\documentclass[titlepage,draft,12pt]{article} 
\usepackage{amsfonts,latexsym} 
\usepackage[a4paper]{geometry}

\date{}

\newcommand{\ep}{\varepsilon}
\newcommand{\qed}{{\penalty 10000\mbox{$\quad\Box$}}}
\newcommand{\re}{\mathbb{R}}

\newcommand{\n}{\mathbb{N}}
\newcommand{\cep}{c_{\ep}}
\newcommand{\gep}{g_{\ep}}
\newcommand{\rep}{r_{\ep}}
\newcommand{\roep}{\rho_{\ep}}
\newcommand{\uep}{u_{\ep}}
\newcommand{\tetep}{\theta_{\ep}}
\newcommand{\da}{D(A)}
\newcommand{\dau}{D(A^{1/2})}
\newcommand{\phibp}{\Phi_{\beta,p}}
\newcommand{\psiap}{\Psi_{\alpha,p}}
\newcommand{\psigp}{\Psi_{\gamma,p}}
\newcommand{\Eep}{E_{\ep}}
\newcommand{\Fep}{F_{\ep}}
\newcommand{\Gep}{G_{\ep}}
\newcommand{\Hep}{H_{\ep}}
\newcommand{\EEep}{\mathcal{E}_{\ep}}
\newcommand{\FFep}{\mathcal{F}_{\ep}}
\newcommand{\GGep}{\mathcal{G}_{\ep}}
\newcommand{\yep}{y_{\ep}}
\newcommand{\zep}{z_{\ep}}
\newcommand{\wep}{w_{\ep}}
\newcommand{\vep}{v_{\ep}}
\newcommand{\lep}{\lambda_{\ep}}
\def\dfrac#1#2{\displaystyle{\frac{#1}{#2}}}

\newtheorem{thm}{Theorem}[section]
\newtheorem{thmbibl}{Theorem}

\newtheorem{rmk}[thm]{Remark}
\newtheorem{prop}[thm]{Proposition}

\newtheorem{lemma}[thm]{Lemma}
\newtheorem{open}[thm]{Open problem}

\title{On the parabolic regime of a hyperbolic equation with weak 
dissipation: the coercive case}

\author{Marina Ghisi\vspace{1ex}\\ 
{\normalsize Universit\`a degli Studi di Pisa} \\
{\normalsize Dipartimento di Matematica ``Leonida Tonelli''}\\ 
{\normalsize PISA (Italy)}\\
{\normalsize e-mail: \texttt{ghisi@dm.unipi.it}}
\and
Massimo Gobbino\vspace{1ex}\\ 
{\normalsize Universit\`a degli Studi di Pisa} \\
{\normalsize Dipartimento di Matematica Applicata ``Ulisse Dini''}\\ 
{\normalsize PISA (Italy)}\\  
{\normalsize e-mail: \texttt{m.gobbino@dma.unipi.it}}}

\begin{document}
\maketitle
\begin{abstract}
	We consider a family of Kirchhoff equations with a small parameter
	$\ep$ in front of the second-order time-derivative, and a
	dissipation term with a coefficient which tends to 0 as $t\to
	+\infty$.
	
	It is well-known that, when the decay of the coefficient is slow 
	enough, solutions behave as solutions of the corresponding 
	parabolic equation, and in particular they decay to 0 as $t\to 
	+\infty$.
	
	In this paper we consider the nondegenerate and coercive case, and
	we prove \emph{optimal} decay estimates for the hyperbolic
	problem, and optimal decay-error estimates for the difference
	between solutions of the hyperbolic and the parabolic problem.
	These estimates show a quite surprising fact: in the coercive case
	the analogy between parabolic equations and dissipative hyperbolic
	equations is weaker than in the noncoercive case.
	
	This is actually a result for the corresponding linear equations
	with time-dependent coefficients.  The nonlinear term comes into
	play only in the last step of the proof.
	
\vspace{1cm}

\noindent{\bf Mathematics Subject Classification 2010 (MSC2010):}
35B25, 35L72, 35B40.


\vspace{1cm} 

\end{abstract}

 
\section{Introduction}

Let $H$ be a real Hilbert space.  For every $x$ and $y$ in
$H$, $|x|$ denotes the norm of $x$, and $\langle x,y\rangle$ denotes
the scalar product of $x$ and $y$.  Let $A$ be a self-adjoint linear
operator on $H$ with dense domain $D(A)$.  We assume that $A$ is
nonnegative, namely $\langle Ax,x\rangle\geq 0$ for every $x\in D(A)$,
so that for every $\alpha\geq 0$ the power $A^{\alpha}x$ is defined
provided that $x$ lies in a suitable domain $D(A^{\alpha})$.

We consider the Cauchy problem
\begin{equation}
	\ep\uep''(t)+\frac{1}{(1+t)^{p}}\uep'(t)+
	m\left(|A^{1/2}\uep(t)|^{2}\right)A\uep(t)=0
	\quad\quad
	\forall t\geq 0,
	\label{pbm:h-eq}
\end{equation}
\begin{equation}
	\uep(0)=u_0,\hspace{3em}\uep'(0)=u_1,
	\label{pbm:h-data}
\end{equation}
where $\ep>0$ and $p\geq 0$ are real parameters,
$m:[0,+\infty)\to(0,+\infty)$ is a locally Lipschitz continuous
function, and $(u_{0},u_{1})\in\da\times\dau$.

The singular perturbation problem in its generality consists in
proving the convergence of solutions of (\ref{pbm:h-eq}),
(\ref{pbm:h-data}) to solutions of the first order problem
\begin{equation}
	\frac{1}{(1+t)^{p}}u'(t)+
	m\left(|A^{1/2}u(t)|^{2}\right)Au(t)=0
	\quad\quad
	\forall t\geq 0,
	\label{pbm:p-eq}
\end{equation}
\begin{equation}
	u(0)=u_{0},
	\label{pbm:p-data}
\end{equation}
obtained setting formally $\ep=0$ in (\ref{pbm:h-eq}), and
omitting the second initial condition in~(\ref{pbm:h-data}).

Several cases have been considered in the last 30 years, depending on
the nonlinearity (degenerate or nondegenerate), on the dissipative
term (constant dissipation $p=0$ or weak dissipation $p>0$), and on
the operator $A$ (coercive or noncoercive).  The main research lines
concern global existence for the parabolic and the hyperbolic problem
(at least when $\ep$ is small enough), decay estimates on $u(t)$,
$\uep(t)$, and $\uep(t)-u(t)$ as $t\to +\infty$, error estimates on
the difference as $\ep\to 0^{+}$, and decay-error estimates, namely
estimates describing in the same time the behavior of the difference
$\uep(t)-u(t)$ as $t\to +\infty$ and $\ep\to 0^{+}$.  The interested
reader is referred to the survey~\cite{gg:survey-diss}, or to the more
recent papers~\cite{ghisi:error,gg:de-dg1,gg:de-dg2}.

In this paper we focus on the case where the equation is
\emph{nondegenerate}, namely
\begin{equation}
	\inf\left\{m(\sigma):\sigma\geq 0\right\}=:\mu>0,
	\label{hp:ndg}
\end{equation}
and  the operator is \emph{coercive}, namely
\begin{equation}
	\inf\left\{\langle Au,u\rangle:u\in D(A),\ |u|=1\right\}=:\nu>0.
	\label{hp:coercive}
\end{equation}

Concerning the parabolic problem, it is well-known that it admits a 
global solution for every $p\geq 0$, and every $u_{0}\in D(A)$ (and 
even for less regular data and more general nonlinearities, 
see~\cite{k-par}).

As for the hyperbolic problem, things are different depending on $p$. 
Let us begin with the linear equation in which $m(\sigma)$ is a positive 
constant. In this case, T.\ Yamazaki~\cite{yamazaki:linear} and J.\ 
Wirth~\cite{wirth} proved two complementary results, which can be 
outlined as follows.
\begin{itemize}
	\item When $p>1$, the dissipative term is too weak, and solutions
	of (\ref{pbm:h-eq}), (\ref{pbm:h-data}) behave as solutions of the
	same equation without the dissipative term.  In particular,
	solutions do not decay to 0.  This is the \emph{hyperbolic
	regime}.

	\item When $p<1$, inertia is negligible, and solutions of
	(\ref{pbm:h-eq}), (\ref{pbm:h-data}) behave as solutions of
	(\ref{pbm:p-eq}), (\ref{pbm:p-data}).  In particular, they decay
	to 0.  This is the \emph{parabolic regime}, with the so-called
	effective dissipation.
	
	\item When $p=1$, the dissipation is still effective (namely the
	integral of the coefficient diverges), but according
	to~\cite{wirth} ``the parabolic asymptotics changes to a wave type
	asymptotics''.  In any case, solutions keep on going to 0, at
	least when $\ep$ is small enough, and for this reason the case
	$p=1$ eventually falls in the parabolic regime.
\end{itemize}

These results have been extended to Kirchhoff equation by H.\
Hashimoto and T.\ Yamazaki~\cite{yamazaki}, T.\
Yamazaki~\cite{yamazaki-wd,yamazaki-cwd} and the
authors~\cite{gg:w-ndg}, in the following sense.
\begin{itemize}
	\item When $p\in[0,1]$, problem (\ref{pbm:h-eq}),
	(\ref{pbm:h-data}) has a unique global solution provided that
	$\ep$ is small enough, and this solution decays to 0 as $t\to
	+\infty$.  This is the parabolic regime.

	\item When $p>1$, existence of global solutions to
	(\ref{pbm:h-eq}), (\ref{pbm:h-data}) is known only for special
	initial data or special operators, the same ones for which global
	existence is known in the nondissipative case.  Global existence
	for every $(u_{0},u_{1})\in D(A)\times D(A^{1/2})$, even for $\ep$
	small enough, is still an open problem, exactly as in the
	nondissipative case.  In any case, nontrivial global solutions, if
	they exist, can \emph{not} decay to 0 as $t\to +\infty$.  This is
	the hyperbolic regime.
\end{itemize}

All the results quoted above do not depend on the coerciveness of $A$,
namely they are true also when $\nu=0$.  

Several estimates on solutions have been proved in the literature,
once again without assumption (\ref{hp:coercive}).  The prototype of
\emph{decay estimates} is that
$$|A^{1/2}u(t)|^{2}\leq\frac{C}{(1+t)^{1+p}}
\quad\mbox{and}\quad
|A^{1/2}\uep(t)|^{2}\leq\frac{C}{(1+t)^{1+p}}$$
for every $t\geq 0$, where the constant $C$ is independent of 
$\ep$ and of course also of $t$. As a consequence, we have also 
that
\begin{equation}
	|A^{1/2}(\uep(t)-u(t))|^{2}\leq\frac{C}{(1+t)^{1+p}}
	\quad\quad
	\forall t\geq 0.
	\label{est:decay-proto}
\end{equation}

The prototype of \emph{error estimates} is that for initial data
$(u_{0},u_{1})\in D(A^{3/2})\times D(A^{1/2})$ one has that
\begin{equation}
	|A^{1/2}(\uep(t)-u(t))|^{2}\leq C\ep^{2} \quad\quad \forall
	t\geq 0,
	\label{est:error-proto}
\end{equation}
where the constant $C$ is once again independent of $\ep$ and $t$
(global-in-time error estimates).  It is well-known that $\ep^{2}$
is the best possible convergence rate (even when looking for
local-in-time error estimates), and that $D(A^{3/2})\times
D(A^{1/2})$ is the minimal requirement on initial data which
guarantees this rate (even in the case of linear equations).  We
refer to~\cite{ch,gg:l-cattaneo,gg:k-PS} for these aspects.

The prototype of \emph{decay-error estimates} is that for initial data
$(u_{0},u_{1})\in D(A^{3/2})\times D(A^{1/2})$ one has that
\begin{equation}
	|A^{1/2}(\uep(t)-u(t))|^{2}\leq 
	C\frac{\ep^{2}}{(1+t)^{1+p}}
	\quad\quad
	\forall t\geq 0.
	\label{est:de-proto}
\end{equation}

We point out in particular that, according to these estimates,
solutions of the hyperbolic problem decay with the same rate of
solutions of the parabolic problem.  Moreover, in the decay-error
estimates (\ref{est:de-proto}) we have the same convergence rate of
the error estimates (\ref{est:error-proto}), and the same decay rate
of the decay estimates (\ref{est:decay-proto}).  Finally, all these
results hold true without coerciveness assumptions on $A$, and for
these general operators it turns out that decay rates are optimal.

When the operator $A$ is coercive, better decay rates are expected.
For example, it is easy to see that solutions of the parabolic problem
satisfy
\begin{equation}
	|A^{1/2}u(t)|^{2}\leq C e^{-\alpha(1+t)^{1+p}}
	\quad\quad
	\forall t\geq 0
	\label{est:decay-p-par}
\end{equation}
for a suitable $\alpha>0$, depending on $\mu$, $\nu$, and $p$ (see
Theorem~\ref{thm:main-p}).

Therefore, the analogy with the noncoercive case could lead to guess
that also solutions of the hyperbolic problem should decay with the
same exponential rate, and the same rate should also appear in the
decay-error estimates.

In this paper we show that this is \emph{not} the case, because
solutions of the hyperbolic problem decay to 0 with a different,
slower rate. Indeed we prove (see Theorem~\ref{thm:main-h}) that
\begin{equation}
	|A^{1/2}\uep(t)|^{2}\leq C e^{-\alpha(1+t)^{1-p}}
	\quad\quad
	\forall t\geq 0
	\label{est:decay-p-hyp}
\end{equation}
if $p\in[0,1)$, and
$$|A^{1/2}\uep(t)|^{2}\leq \frac{C}{(1+t)^{\alpha}}
\quad\quad
\forall t\geq 0
$$
if $p=1$, where $\alpha<2\mu\nu$ if $p=0$, and $\alpha$ is any
(positive) real number if $p\in(0,1]$ (now the constant $C$ depends
also on $\alpha$).  These rates are optimal, in the sense that every
nonzero solution does not satisfy an estimate such as
(\ref{est:decay-p-hyp}) with an exponent larger than $(1-p)$ (see
Theorem~\ref{thm:main-opt}).  The same slower rates appear also in the
decay-error estimates (see Theorem~\ref{thm:main-sp}), and of course
they are optimal also in this case.

We have thus shown an essential difference between the coercive and
the noncoercive case.  In the noncoercive case, solutions of the
hyperbolic problem mimic the behavior of solutions of the parabolic
problem for every $p\in[0,1]$.  In the coercive case, this is true
only for $p=0$, when the exponent $(1+p)$ in (\ref{est:decay-p-par})
and the exponent $(1-p)$ in (\ref{est:decay-p-hyp}) coincide.  On the
contrary, for every $p\in(0,1]$ there is a spread between exponents in
the decay rates of $u(t)$ and $\uep(t)$, and this spread becomes
larger and larger as $p$ approaches 1.  As a consequence, from the
point of view of decay rates, (\ref{pbm:p-eq}) is a good approximation
of (\ref{pbm:h-eq}) for $\ep$ small in the noncoercive case, but not
in the coercive case (see also section~\ref{sec:heuristics}).

In both cases (coercive and noncoercive), the parabolic problem and
the hyperbolic problem take different paths when $p>1$: solutions of
the parabolic problem keep on decaying according to
(\ref{est:decay-p-par}), hence faster and faster as $p$ grows, while
solutions of the hyperbolic problem do not decay to 0 any more
(provided that they globally exist).

All our proofs are based on linear arguments.  To this end, we
first linearize (\ref{pbm:h-eq}) and (\ref{pbm:p-eq}).  We obtain the
following equations
\begin{equation}
	\ep\uep''(t)+\frac{1}{(1+t)^{p}}\uep'(t)+
	\cep(t)A\uep(t)=0
	\quad\quad
	\forall t\geq 0,
	\label{pbm:h-leq}
\end{equation}
\begin{equation}
	\frac{1}{(1+t)^{p}}u'(t)+c(t)Au(t)=0
	\quad\quad
	\forall t\geq 0,
	\label{pbm:p-leq}
\end{equation}
with time-dependent coefficients $\cep:[0,+\infty)\to(0,+\infty)$ and 
$c:[0,+\infty)\to(0,+\infty)$. 

Then we prove decay and decay-error estimates for solutions of these
linear equations, under suitable assumptions on the coefficients.
This is the core of the paper.

Finally, we just observe that the coefficients $\cep(t)$ and $c(t)$
coming from the nonlinear terms in (\ref{pbm:h-eq}) and
(\ref{pbm:p-eq}) satisfy the assumptions required by the linear
theory.  Fortunately, these assumptions are quite weak, and follow
easily from previous literature on the noncoercive case.

This paper is organized as follows.  In section~\ref{sec:bibl} we
recall the previous results and estimates needed throughout this
paper.  In section~\ref{sec:main} we state our main results for
Kirchhoff equations.  In section~\ref{sec:heuristics} we present a
heuristic argument leading to our decay rates.  In
section~\ref{sec:linear} we state our results for linear equations
with time-dependent coefficients.  In section~\ref{sec:proofs} we
collect all proofs.

\setcounter{equation}{0}
\section{Statements}\label{sec:statements}

\subsection{Previous works}\label{sec:bibl}

The theory of nondegenerate Kirchhoff equations with weak dissipation
has been developed in \cite{yamazaki-wd,yamazaki-cwd,gg:w-ndg}.  In
the following statement we collect the existence results, and some
decay and error estimates.  We limit ourselves to the results which
are needed in the sequel, and for this reason Theorem~\ref{thm:bibl}
below does not represent the full state of the art, especially for
decay-error estimates.  The interested reader is referred to section~5
of~\cite{gg:survey-diss} for further (and more refined) estimates and
references.

\begin{thmbibl}\label{thm:bibl}
	Let $H$ be a Hilbert space, let $A$ be a self-adjoint nonnegative
	operator on $H$ with dense domain $D(A)$ (no coercivity assumption
	on $A$), let $m:[0,+\infty)\to(0,+\infty)$ be a locally Lipschitz
	continuous function satisfying the nondegeneracy condition
	(\ref{hp:ndg}), and let $(u_{0},u_{1})\in D(A)\times D(A^{1/2})$.
	
	Then we have the following conclusions.
	\begin{enumerate}
		\renewcommand{\labelenumi}{(\arabic{enumi})} 
		\item \emph{(Parabolic problem)} For every $p\geq 0$, problem
		(\ref{pbm:p-eq}), (\ref{pbm:p-data}) has a unique global
		solution
		\begin{equation}
			u\in C^{1}\left([0,+\infty);H\right)\cap
			C^{0}\left([0,+\infty);\da\right).
			\label{th:p-reg}
		\end{equation}
	
		Moreover $u\in C^{1}\left((0,+\infty);D(A^{\alpha})\right)$
		for every $\alpha\geq 0$ (and more generally $u$ is of class
		$C^{k+1}$ when $m(\sigma)$ is of class $C^{k}$), and there
		exists a constant $C$ such that
		\begin{equation}
			(1+t)^{2}|u'(t)|^{2}+(1+t)^{1+p}|A^{1/2}u(t)|^{2}+
			(1+t)^{2(1+p)}|Au(t)|^{2}\leq C
			\quad\quad
			\forall t\geq 0.
			\label{th:bibl-p}
		\end{equation}
	
		\item \emph{(Hyperbolic problem)} For every $p\in[0,1]$, there
		exists $\ep_{0}>0$ such that, for every $\ep\in(0,\ep_{0})$,
		problem (\ref{pbm:h-eq}), (\ref{pbm:h-data}) has a unique
		global solution
		\begin{equation}
			\uep\in C^{2}\left([0,+\infty);H\right)\cap
			C^{1}\left([0,+\infty);\dau\right)\cap
			C^{0}\left([0,+\infty);\da\right).
			\label{th:h-reg}
		\end{equation}
		
		Moreover, there exists a constant $C$ such that for every 
		$\ep\in(0,\ep_{0})$ we have that
		\begin{equation}
			(1+t)^{2}|\uep'(t)|^{2}+(1+t)^{1+p}|A^{1/2}\uep(t)|^{2}+
			(1+t)^{2(1+p)}|A\uep(t)|^{2}\leq C
			\quad\quad
			\forall t\geq 0.
			\label{th:bibl-h}
		\end{equation}
	
		\item \emph{(Singular perturbation)} If $p\in[0,1]$, and
		$(u_{0},u_{1})\in D(A^{3/2})\times D(A^{1/2})$, then there
		exist $\ep_{1}\in(0,\ep_{0})$ and $C$ such that, for every
		$\ep\in(0,\ep_{1})$ we have that
		\begin{equation}
			|A^{1/2}(\uep(t)-u(t))|^{2}\leq C\ep^{2}
			\quad\quad
			\forall t\geq 0.
			\label{th:bibl-sp}
		\end{equation}
	\end{enumerate}
\end{thmbibl}

\subsection{Main results}\label{sec:main}

In this section we state the main results of this paper. The first 
one concerns decay estimates for solutions of the parabolic problem.

\begin{thm}[Parabolic equation]\label{thm:main-p}
	Let $H$ be a Hilbert space, and let $A$ be a self-adjoint operator
	on $H$ with dense domain $D(A)$.  Let $u_{0}\in D(A)$, let $p\geq
	0$, and let $m:[0,+\infty)\to(0,+\infty)$ be a locally Lipschitz
	continuous function.
	
	Let us assume that the nondegeneracy and coerciveness assumptions
	(\ref{hp:ndg}) and (\ref{hp:coercive}) are satisfied.
	
	Then problem (\ref{pbm:p-eq}), (\ref{pbm:p-data}) has a unique
	global solution $u(t)$ with the regularity prescribed in
	statement~(1) of Theorem~\ref{thm:bibl}, and there exists a
	constant $C$ such that
	\begin{equation}
		|u(t)|^{2}+|A^{1/2}u(t)|^{2}+|Au(t)|^{2}+
		\frac{|u'(t)|^{2}}{(1+t)^{2p}}\leq
		C\exp\left(-\frac{2\mu\nu}{1+p}(1+t)^{1+p}\right)
		\label{th:p}
	\end{equation}
	for every $t\geq 0$.
\end{thm}

The second result concerns decay estimates for solutions of the
hyperbolic problem.

\begin{thm}[Hyperbolic equation]\label{thm:main-h}
	Let $H$ be a Hilbert space, and let $A$ be a self-adjoint operator
	on $H$ with dense domain $D(A)$.  Let $(u_{0},u_{1})\in D(A)\times
	D(A^{1/2})$, let $p\in[0,1]$, and let
	$m:[0,+\infty)\to(0,+\infty)$ be a locally Lipschitz continuous
	function.
	
	Let us assume that the nondegeneracy and coerciveness assumptions
	(\ref{hp:ndg}) and (\ref{hp:coercive}) are satisfied.
	
	Then there exists $\ep_{0}>0$ such that, for every
	$\ep\in(0,\ep_{0})$, problem (\ref{pbm:h-eq}), (\ref{pbm:h-data})
	has a unique global solution $\uep(t)$ with the regularity 
	prescribed by (\ref{th:h-reg}).
	
	Moreover the function
	\begin{equation}
		\Gamma_{\ep}(t):=|\uep(t)|^{2}+|A^{1/2}\uep(t)|^{2}+|A\uep(t)|^{2}+
		|\uep'(t)|^{2}+\ep|A^{1/2}\uep'(t)|^{2}
		\label{defn:Gamma-ep}
	\end{equation}
	satisfies the following decay estimates.
	\begin{itemize}
		\item \fbox{Case $p=0$} For every $\beta<2\mu\nu$, there exist
		$\ep_{1}\in(0,\ep_{0}]$ and $C$ such that
		\begin{equation}
			\Gamma_{\ep}(t)\leq
			Ce^{-\beta t}
			\quad\quad
			\forall t\geq 0,\ \forall\ep\in(0,\ep_{1}).
			\label{th:h-p=0}
		\end{equation}
	
		\item  \fbox{Case $p\in(0,1)$} For every $\beta>0$, there exist 
		$\ep_{1}\in(0,\ep_{0}]$ and $C$ such that 
		\begin{equation}
			\Gamma_{\ep}(t)\leq
			Ce^{-\beta(1+t)^{1-p}}
			\quad\quad
			\forall t\geq 0,\ \forall\ep\in(0,\ep_{1}).
			\label{th:h-p=(0,1)}
		\end{equation}
	
		\item  \fbox{Case $p=1$} For every $\beta>0$, there exist 
		$\ep_{1}\in(0,\ep_{0}]$ and $C$ such that 
		\begin{equation}
			\Gamma_{\ep}(t)\leq
			\frac{C}{(1+t)^{\beta}}
			\quad\quad
			\forall t\geq 0,\ \forall\ep\in(0,\ep_{1}).
			\label{th:h-p=1}
		\end{equation}
	\end{itemize}
\end{thm}

Of course the constants $C$ and $\ep_{1}$ in (\ref{th:h-p=0}) through
(\ref{th:h-p=1}) depend also on $\beta$.

The third step concerns the singular perturbation problem.  Following
the approach introduced in~\cite{lions} in the linear case, we define
the corrector $\tetep(t)$ as the solution of the second order
\emph{linear} ordinary differential equation
\begin{equation}
	\ep\tetep''(t)+\frac{1}{(1+t)^{p}}\tetep'(t)=0 \hspace{2em}
	\forall t\geq 0,
	\label{pbm:tetep-eq}
\end{equation}
with initial data
$$\tetep(0)=0,\hspace{2em}\tetep'(0)=u_1+
m\left(|A^{1/2}u_{0}|^{2}\right)Au_{0}=:\theta_{0}.$$

Since $\theta_{0}=\uep'(0)-u'(0)$, this corrector keeps into account
the boundary layer due to the loss of one initial condition.  

We can now define $\rep(t)$ and $\roep(t)$ in such a way that
$$\uep(t)=u(t)+\tetep(t)+\rep(t)=u(t)+\roep(t)\quad\quad\forall t\geq
0.$$

With these notations, the singular perturbation problem consists in
proving that $\rep(t)\to 0$ or $\roep(t)\to 0$ in some sense as
$\ep\to 0^{+}$.  We recall that the two remainders play different
roles.  In particular, $\rep(t)$ is well suited for estimating
derivatives, while $\roep(t)$ is used in estimates without
derivatives.  This distinction is essential.  Indeed it is not
possible to prove decay-error estimates on $A^{\alpha}\rep(t)$ because
it does not decay to 0 as $t\to +\infty$ (indeed $\uep(t)$ and $u(t)$
tend to 0, while the corrector $\tetep(t)$ does not), and it is not
possible to prove decay-error estimates on $A^{\alpha}\roep'(t)$
because in general for $t=0$ it does not tend to 0 as $\ep\to 0^{+}$
(due to the loss of one initial condition).

We are now ready to state our decay-error estimates.

\begin{thm}[Singular perturbation]\label{thm:main-sp}
	Let $H$ be a Hilbert space, and let $A$ be a self-adjoint operator
	on $H$ with dense domain $D(A)$.  Let $(u_{0},u_{1})\in D(A)\times
	D(A^{1/2})$, let $p\in[0,1]$, and let
	$m:[0,+\infty)\to(0,+\infty)$ be a locally Lipschitz continuous
	function.
	
	Let us assume that the nondegeneracy and coerciveness assumptions
	(\ref{hp:ndg}) and (\ref{hp:coercive}) are satisfied, and let
	$u(t)$, $\ep_{0}$, $\uep(t)$, $\rep(t)$, $\roep(t)$ be as above.
	
	Let us consider the functions
	$$\Gamma_{r,\ep}(t):=|\roep(t)|^{2}+|A^{1/2}\roep(t)|^{2}+
	\ep|\rep'(t)|^{2},$$
	$$\Gamma_{c,\ep}(t):=|\roep(t)|^{2}+|A^{1/2}\roep(t)|^{2}+|A\roep(t)|^{2}+
	|\rep'(t)|^{2}+\ep|A^{1/2}\rep'(t)|^{2},$$
	where indices $c$ and $r$ stay for ``complete'', and ``reduced'', 
	respectively.
	\begin{enumerate}
		\renewcommand{\labelenumi}{(\arabic{enumi})}
		\item  If in addition $(u_{0},u_{1})\in D(A^{3/2})\times 
		D(A^{1/2})$, then we have the following decay-error estimates.
		\begin{itemize}
			\item  \fbox{Case $p=0$} For every $\beta<2\mu\nu$, there exist 
			$\ep_{1}\in(0,\ep_{0}]$ and $C$ such that
			\begin{equation}
				\Gamma_{r,\ep}(t)\leq
				C\ep^{2}e^{-\beta t}
				\quad\quad
				\forall t\geq 0,\ \forall\ep\in(0,\ep_{1}).
				\label{th:sp-p=0}
			\end{equation}
		
			\item  \fbox{Case $p\in(0,1)$} For every $\beta>0$, there exist 
			$\ep_{1}\in(0,\ep_{0}]$ and $C$ such that
			\begin{equation}
				\Gamma_{r,\ep}(t)\leq
				C\ep^{2}e^{-\beta(1+t)^{1-p}}
				\quad\quad
				\forall t\geq 0,\ \forall\ep\in(0,\ep_{1}).
				\label{th:sp-p=0-1}
			\end{equation}
		
			\item  \fbox{Case $p=1$} For every $\beta>0$, there exist 
			$\ep_{1}\in(0,\ep_{0}]$ and $C$ such that
			\begin{equation}
				\Gamma_{r,\ep}(t)\leq
				\frac{C\ep^{2}}{(1+t)^{\beta}}
				\quad\quad
				\forall t\geq 0,\ \forall\ep\in(0,\ep_{1}).
				\label{th:sp-p=1}
			\end{equation}
		\end{itemize}
	
		\item  If in addition $(u_{0},u_{1})\in D(A^{2})\times 
		D(A)$, then we have the same decay-error estimates with 
		$\Gamma_{c,\ep}(t)$ instead of $\Gamma_{r,\ep}(t)$.
	\end{enumerate}
\end{thm}

As in Theorem~\ref{thm:main-h} above, the constants $C$ and $\ep_{1}$
in (\ref{th:sp-p=0}) through (\ref{th:sp-p=1}) depend also on
$\beta$.  We point out that in these estimates we have the same
convergence rate as in (\ref{th:bibl-sp}), and the same decay rates as
in (\ref{th:h-p=0}) through (\ref{th:h-p=1}).

The last result we state, together with Remarks~\ref{rmk:opt-p>0}
and~\ref{rmk:opt-p=0} below, clarifies the optimality of the decay
rates of Theorem~\ref{thm:main-h}, hence also of
Theorem~\ref{thm:main-sp}.

\begin{thm}[Optimality of decay rates]\label{thm:main-opt}
	Let $H$, $A$, $p\in[0,1]$, $m:[0,+\infty)\to(0,+\infty)$, and
	$(u_{0},u_{1})\in D(A)\times D(A^{1/2})$ be as in
	Theorem~\ref{thm:main-h}.  Let $\ep>0$, and let $\uep(t)$ be the
	solution to problem (\ref{pbm:h-eq}), (\ref{pbm:h-data}).
	
	Let $\Phi:[0,+\infty)\to(0,+\infty)$ be a function of class 
	$C^{1}$ such that
	\begin{equation}
		\lim_{t\to +\infty}(1+t)^{p}\frac{\Phi'(t)}{\Phi(t)}=-\infty.
		\label{hp:main-opt}
	\end{equation}
	
	If $(u_{0},u_{1})\neq(0,0)$, then
	\begin{equation}
		\lim_{t\to +\infty}\left(
		\ep|\uep'(t)|^{2}+|A^{1/2}\uep(t)|^{2}\right)
		\frac{1}{\Phi(t)}=+\infty.
		\label{th:main-opt}
	\end{equation}
\end{thm}

\begin{rmk}\label{rmk:opt-p>0}
	\begin{em}
		When $p>0$, Theorem~\ref{thm:main-opt} is exactly the
		counterpart of Theorem~\ref{thm:main-h}.  Indeed let us
		consider any $\Phi:[0,+\infty)\to[0,+\infty)$, and let
		$\Gamma_{\ep}(t)$ be defined as in (\ref{defn:Gamma-ep}).  If
		(\ref{hp:main-opt}) is satisfied, then we can not expect that
		$\Gamma_{\ep}(t)\leq C\Phi(t)$ because of (\ref{th:main-opt}).
		On the contrary, if
		$$(1+t)^{p}\frac{\Phi'(t)}{\Phi(t)}\geq-\beta>-\infty,$$
		then $\Phi(t)\geq Ce^{-\beta(1-p)^{-1}(1+t)^{1-p}}$ if
		$p\in(0,1)$, and $\Phi(t)\geq C(1+t)^{-\beta}$ if $p=1$, and
		in both cases Theorem~\ref{thm:main-h} guarantees that
		$\Gamma_{\ep}(t)\leq C\Phi(t)$.
		
		Note in particular that the function
		$\Phi(t):=e^{-\beta(1+t)^{\delta}}$ satisfies 
		(\ref{hp:main-opt}) if and only if $\delta>1-p$, which means 
		that $(1-p)$ is the larger exponent for which 
		(\ref{th:h-p=(0,1)}) holds true.
	\end{em}
\end{rmk}

\begin{rmk}\label{rmk:opt-p=0}
	\begin{em}
		When $p=0$, estimate (\ref{th:h-p=0}) can not be true when 
		$\beta>2\mu\nu$. This can be easily seen by considering the 
		explicit solutions of the ordinary differential equation
		\begin{equation}
			\ep y''(t)+y'(t)+\mu\nu y(t)=0,
			\label{pbm:h-ODE}
		\end{equation}
		which is just the particular case of (\ref{pbm:h-eq}) where 
		$H=\re$, $A$ is $\nu$ times the identity, and 
		$m(\sigma)\equiv\mu$ is a constant.
		
		On the other hand, solutions of (\ref{pbm:h-ODE}) satisfy
		(\ref{th:h-p=0}) also with $\beta=2\mu\nu$.  We suspect that
		this could be true in general, but for the time being we have
		no proof.
	\end{em}
\end{rmk}

\begin{open}
	\begin{em}
		Is (\ref{th:h-p=0}) true also in the case $\beta=2\mu\nu$?
	\end{em}
\end{open}

\subsection{Heuristics}\label{sec:heuristics}

According to Theorem~\ref{thm:main-p}, solutions of the
parabolic problem decay as solutions of the ordinary differential
equation
\begin{equation}
	\frac{1}{(1+t)^{p}}y'(t)+\mu\nu\, y(t)=0.
	\label{pbm:toy-p}
\end{equation}

This is hardly surprising, since (\ref{pbm:toy-p}) is just the special
case of (\ref{pbm:p-eq}) corresponding to $H=\re$, $A$ equal to $\nu$
times the identity, and $m(\sigma)\equiv\mu$.

Analogously, it is reasonable to expect solutions of the hyperbolic
problem to decay as solutions of the ordinary differential equation
\begin{equation}
	\ep\yep''(t)+\frac{1}{(1+t)^{p}}\yep'(t)+\mu\nu\, \yep(t)=0.
	\label{pbm:toy-h}
\end{equation}

A reasonable ansatz for these solutions is that asymptotically they 
are the product of an oscillatory term $\vep(t)$, and a decaying term 
$\lep(t)$. Plugging $\yep(t)=\lep(t)\cdot\vep(t)$ into 
(\ref{pbm:toy-h}), we obtain that
$$\left(\strut\ep\vep''(t)+\mu\nu\, \vep(t)\right)\lep(t)+
\left(2\ep\lep'(t)+\frac{\lep(t)}{(1+t)^{p}}\right)\vep'(t)+
\left(\ep\lep''(t)+\frac{\lep'(t)}{(1+t)^{p}}\right)\vep(t) = 0.$$

A reasonable guess is now that the coefficient of $\lep(t)$ in the
first term is almost 0, as well as the coefficient of $\vep'(t)$ in
the second term.

The first condition is that $\ep\vep''(t)+\mu\nu\vep(t)\sim 0$, namely 
$$\vep(t)\sim\sin\left(\sqrt{\frac{\mu\nu}{\ep}}t\right),$$
which yields the same oscillations of the undamped equation.

The second condition is that
\begin{equation}
	2\ep\lep'(t)+\frac{\lep(t)}{(1+t)^{p}}\sim0,
	\label{pbm:toy-lep}
\end{equation}
and for every $p\in(0,1]$ this yields a decay rate which is 
compatible both with Theorem~\ref{thm:main-h} and with 
Theorem~\ref{thm:main-opt}. 

We do not know if similar asymptotics have been rigorously justified
in the literature (see~\cite{wirth} for the case $p=1$).
Nevertheless, this non-rigorous argument suggests that actually there
is no sharp break between parabolic and hyperbolic regimes.  For
$p\leq 1$, the hyperbolic nature survives in the oscillatory behavior
of $\vep(t)$, but it is hidden by the damping imposed by
(\ref{pbm:toy-lep}).  When $p>1$, solutions of (\ref{pbm:toy-lep})
tend to a positive constant, and the hyperbolic nature emerges
undisputed.

We conclude by pointing out once again that this analysis applies to
the nondegenerate coercive case.  Things are quite different both in
the nondegenerate noncoercive case
(see~\cite{gg:w-ndg,wirth-jde,wirth,yamazaki:linear,yamazaki-wd}), and
in the degenerate coercive case (see~\cite{gg:de-dg1,gg:de-dg2}).

\subsection{Linearization}\label{sec:linear}

Proofs of our main results are based on the analysis of the linear 
equations (\ref{pbm:h-leq}) and (\ref{pbm:p-leq}). We assume that 
the coefficient $c:[0,+\infty)\to(0,+\infty)$ is of class 
$C^{1}$, and satisfies the following estimates
\begin{equation}
	c(t)\geq \mu>0
	\quad\quad
	\forall t\geq 0,
	\label{hp:c-below}
\end{equation}
\begin{equation}
	c(t)\leq M_{1}
	\quad\quad
	\forall t\geq 0,
	\label{hp:c-above}
\end{equation}
\begin{equation}
	|c'(t)|\leq M_{2}
	\quad\quad
	\forall t\geq 0.
	\label{hp:c'-above}
\end{equation}

Similarly, we assume that $\cep:[0,+\infty)\to(0,+\infty)$, with
$\ep\in(0,\ep_{0})$, is a family of coefficients of class $C^{1}$
satisfying the following estimates
\begin{equation}
	\cep(t)\geq\mu>0
	\quad\quad
	\forall t\geq 0,\ \forall\ep\in(0,\ep_{0}),
	\label{hp:cep-below}
\end{equation}
\begin{equation}
	\cep(t)\leq M_{3}
	\quad\quad
	\forall t\geq 0,\ \forall\ep\in(0,\ep_{0}),
	\label{hp:cep-above}
\end{equation}
\begin{equation}
	|\cep'(t)|\leq \frac{M_{4}}{(1+t)^{p}}
	\quad\quad
	\forall t\geq 0,\ \forall\ep\in(0,\ep_{0}).
	\label{hp:cep'-above}
\end{equation}

When considering the singular perturbation, we also assume that
\begin{equation}
	|\cep(t)-c(t)|\leq M_{5}\ep
	\quad\quad
	\forall t\geq 0,\ \forall\ep\in(0,\ep_{0}),
	\label{hp:cep-c}
\end{equation}
and we define the corrector $\tetep(t)$ as the solution of 
(\ref{pbm:tetep-eq}) with initial data
\begin{equation}
	\tetep(0)=0,\hspace{2em}\tetep'(0)=u_1+c(0)Au_{0}=:\theta_{0}.
	\label{pbm:tetep-ldata}
\end{equation}

The following results are the linear counterparts of
Theorems~\ref{thm:main-p} through \ref{thm:main-opt}.  All of them can
be extended to Lipschitz continuous coefficients through a
straightforward approximation argument.

\begin{thm}[Linear parabolic equation]\label{thm:lin-p}
	Let $H$, $A$, $p\geq 0$, and $u_{0}\in D(A)$ be as in
	Theorem~\ref{thm:main-p}.  Let $c:[0,+\infty)\to(0,+\infty)$ be a
	continuous function satisfying (\ref{hp:c-below}) and
	(\ref{hp:c-above}).
	
	Then problem (\ref{pbm:p-leq}), (\ref{pbm:p-data}) has a unique
	global solution $u(t)$ with the regularity prescribed by 
	(\ref{th:p-reg}), and this solution satisfies (\ref{th:p}).

\end{thm}

\begin{thm}[Linear hyperbolic equation]\label{thm:lin-h}
	Let $H$, $A$, $p\in[0,1]$, and $(u_{0},u_{1})\in D(A)\times
	D(A^{1/2})$ be as in Theorem~\ref{thm:main-h}, and let
	$\ep_{0}>0$.  Let $\cep:[0,+\infty)\to(0,+\infty)$, with
	$\ep\in(0,\ep_{0})$, be a family of coefficients of class $C^{1}$
	satisfying (\ref{hp:cep-below}) through (\ref{hp:cep'-above}).
	
	Then, for every $\ep\in(0,\ep_{0})$, problem (\ref{pbm:h-leq}),
	(\ref{pbm:h-data}) has a unique global solution $\uep(t)$ with the
	regularity prescribed by (\ref{th:h-reg}), and this solution
	satisfies the same decay estimates stated in
	Theorem~\ref{thm:main-h}, depending on the values of $p$.
\end{thm}

\begin{thm}[Linear singular perturbation]\label{thm:lin-sp}
	Let $H$, $A$, $p\in[0,1]$, $(u_{0},u_{1})$, $\ep_{0}$, $c(t)$, 
	$u(t)$, $\cep(t)$, $\uep(t)$ be as in Theorems~\ref{thm:lin-p} 
	and~\ref{thm:lin-h}.  
	
	Let us assume that also (\ref{hp:c'-above}) and (\ref{hp:cep-c})
	hold true, and let $\rep(t)$ and $\roep(t)$ be defined as usual
	(keeping in mind that the corrector now satisfies 
	(\ref{pbm:tetep-eq}) and (\ref{pbm:tetep-ldata})).
	
	Then $\rep(t)$ and $\roep(t)$ satisfy the decay-error estimates of
	statements~(1) and~(2) of Theorem~\ref{thm:main-sp}, depending on
	the further regularity of $(u_{0},u_{1})$, and on the values of
	$p$.
\end{thm}

\begin{thm}[Linear hyperbolic equation: optimality]\label{thm:lin-opt}
	Let $H$, $A$, $p\in[0,1]$, and $(u_{0},u_{1})\in D(A)\times
	D(A^{1/2})$ be as in Theorem~\ref{thm:main-h}.  Let $\ep>0$, and
	let $\uep(t)$ be the solution to problem (\ref{pbm:h-leq}),
	(\ref{pbm:h-data}) with a coefficient
	$\cep:[0,+\infty)\to(0,+\infty)$ of class $C^{1}$ satisfying
	(\ref{hp:cep-below}) through (\ref{hp:cep'-above}).
	
	If $(u_{0},u_{1})\neq(0,0)$, then (\ref{th:main-opt}) holds true
	for every function $\Phi:[0,+\infty)\to(0,+\infty)$ of class
	$C^{1}$ satisfying (\ref{hp:main-opt}).
\end{thm}

\setcounter{equation}{0}
\section{Proofs}\label{sec:proofs}

\subsection{Proof of Theorem~\ref{thm:lin-p}}

We prove a more general result, with some further estimates needed 
when dealing with the singular perturbation problem.

\begin{prop}\label{prop:p-lin}
	Let $H$, $A$, and $c(t)$ be as in
	Theorem~\ref{thm:lin-p}.  Let us set
	\begin{equation}
		\gamma:=\frac{2\mu\nu}{1+p},
		\hspace{3em}
		\psiap(t):=\exp\left(
		-\alpha\left[(1+t)^{1+p}-1\right]\right).		
		\label{defn:psiap}
	\end{equation}
	
	Then we have the following estimates.
	\begin{enumerate}
		\renewcommand{\labelenumi}{(\arabic{enumi})}
		\item  If $u_{0}\in D(A^{k/2})$ for some $k\in\n$, then
		\begin{equation}
			|A^{k/2}u(t)|^{2}\leq |A^{k/2}u_{0}|^{2}\psigp(t)
			\quad\quad
			\forall t\geq 0.
			\label{th:p-ak}
		\end{equation}
		
		Moreover, for every $\alpha<\gamma$ we have that
		\begin{equation}
			\int_{0}^{+\infty}\frac{|A^{(k+1)/2}u(t)|^{2}}{\psiap(t)}\,dt
			\leq\left(2\mu-\frac{\alpha(1+p)}{\nu}\right)^{-1}
			|A^{k/2}u_{0}|^{2}.
			\label{th:p-ak-int}
		\end{equation}
	
		\item If $u_{0}\in D(A^{3/2})$, and $c(t)$ is of class $C^{1}$
		and satisfies (\ref{hp:c'-above}), then for every
		$\alpha<\gamma$ there exists a constant $C$ (depending also on
		$\alpha$) such that
		\begin{equation}
			\int_{0}^{+\infty}\frac{|u''(t)|^{2}}{\psiap(t)}\,dt
			\leq C.
			\label{th:u''-int-1}
		\end{equation}
	
		\item If $u_{0}\in D(A^{2})$, and $c(t)$ is of class $C^{1}$
		and satisfies (\ref{hp:c'-above}), then there exists a
		constant $C$ such that
		\begin{equation}
			|u''(t)|^{2}\leq C(1+t)^{4p}\psigp(t)
			\quad\quad
			\forall t\geq 0.
			\label{th:u''-pointwise}
		\end{equation}
		
		Moreover, for every $\alpha<\gamma$, there exists a constant
		$C$ (depending also on $\alpha$) such that
		\begin{equation}
			\int_{0}^{+\infty}\frac{|A^{1/2}u''(t)|^{2}}{\psiap(t)}\,dt
			\leq C.
			\label{th:u''-int-2}
		\end{equation}
	\end{enumerate}
\end{prop}

\paragraph{\textit{\textmd{Proof}}}

Let us set $E_{k}(t):=|A^{k/2}u(t)|^{2}$. From (\ref{pbm:p-leq}), 
(\ref{hp:coercive}), and (\ref{hp:c-below}), we have that
$$E_{k}'(t)=2\langle A^{(k+1)/2}u(t),A^{(k-1)/2}u'(t)\rangle=
-2c(t)(1+t)^{p}\left|A^{(k+1)/2}u(t)\right|^{2}$$
$$\hspace{3em}\leq
-2 c(t)(1+t)^{p}\cdot\nu\left|A^{k/2}u(t)\right|^{2}\leq
-2\mu\nu(1+t)^{p}E_{k}(t).$$

Integrating this differential inequality, we obtain (\ref{th:p-ak}).

Moreover we have that
\begin{eqnarray*}
	\frac{d}{dt}\left[\frac{E_{k}(t)}{\psiap(t)}\right] & = &
	\frac{E_{k}'(t)}{\psiap(t)}
	+\alpha(1+p)(1+t)^{p}\frac{|A^{k/2}u(t)|^{2}}{\psiap(t)} \\
	 & \leq &
	-2\mu(1+t)^{p}\frac{|A^{(k+1)/2}u(t)|^{2}}{\psiap(t)}
	+\frac{\alpha(1+p)}{\nu}(1+t)^{p}
	\frac{|A^{(k+1)/2}u(t)|^{2}}{\psiap(t)},
\end{eqnarray*}
hence
$$\left(2\mu-\frac{\alpha(1+p)}{\nu}\right)
\int_{0}^{t}(1+s)^{p}\frac{|A^{(k+1)/2}u(s)|^{2}}{\psiap(s)}\,ds
+\frac{E_{k}(t)}{\psiap(t)}\leq E_{k}(0)
\quad\quad
\forall t\geq 0,$$
which easily implies (\ref{th:p-ak-int}).

Let us prove the estimates on the second derivative. From 
(\ref{pbm:p-leq}) we obtain that
$$u''(t)=-p(1+t)^{p-1}c(t)Au(t)-(1+t)^{p}c'(t)Au(t)
+(1+t)^{2p}c^{2}(t)A^{2}u(t).$$

Therefore, from (\ref{hp:c-above}) and (\ref{hp:c'-above}), it 
follows that
\begin{equation}
	|u''(t)|^{2}\leq k_{1}(1+t)^{2p}|Au(t)|^{2}
	+k_{2}(1+t)^{4p}|A^{2}u(t)|^{2}.
	\label{comp:u''}
\end{equation}

If $u_{0}\in D(A^{2})$, then (\ref{th:u''-pointwise}) follows from 
(\ref{th:p-ak}) with $k=2$ and $k=4$.

In order to prove the integral estimates on $u''(t)$, let us choose
$\eta$ such that $\alpha<\alpha+\eta<\gamma$. Since
$\Psi_{\alpha+\eta,p}(t)=\psiap(t)\cdot\Psi_{\eta,p}(t)$, and since
$$\sup_{t\geq 0}\left\{\Psi_{\eta,p}(t)(1+t)^{4p}\right\}<+\infty,$$
from (\ref{comp:u''}) it follows that
$$\frac{|u''(t)|^{2}}{\psiap(t)}\leq (1+t)^{4p}\Psi_{\eta,p}(t)\cdot
\frac{k_{1}|Au(t)|^{2} +k_{2}|A^{2}u(t)|^{2}}{\psiap(t)\cdot\Psi_{\eta,p}(t)} \leq
k_{3}\frac{|Au(t)|^{2}+|A^{2}u(t)|^{2}}{\Psi_{\alpha+\eta,p}(t)}.$$

From (\ref{th:p-ak-int}) with $k=1$ and $k=3$ we conclude that
$$\int_{0}^{+\infty}\frac{|u''(t)|^{2}}{\psiap(t)}\,dt\leq
k_{3}\int_{0}^{+\infty}
\frac{|Au(t)|^{2}+|A^{2}u(t)|^{2}}{\Psi_{\alpha+\eta,p}(t)}\,dt
\leq k_{4}$$
for a suitable $k_{4}$ depending also on $\eta$. This proves
(\ref{th:u''-int-1}).

The proof of (\ref{th:u''-int-2}) is completely analogous.\qed

\subsection{Comparison results for ODEs}

In this subsection we prove estimates for solutions of three ordinary
differential equations needed in the sequel.  To begin with, for every
$\beta>0$ and every $p\geq 0$ we define
$\phibp:[0,+\infty)\to(0,+\infty)$ as the solution of the Cauchy
problem
\begin{equation}
	\phibp'(t)=-\frac{\beta}{(1+t)^{p}}\phibp(t)
	\quad\quad
	\forall t\geq 0,
	\label{defn:phibp-eq}
\end{equation}
\begin{equation}
	\phibp(0)=1.
	\label{defn:phibp-0}
\end{equation}

We point out that solutions of this problem decay as the
right-hand sides of (\ref{th:h-p=0}) through (\ref{th:h-p=1}),
depending on the values of $p$.  This is the reason why we are going
to exploit $\phibp(t)$ several times in the proofs of our decay and
decay-error estimates.

\begin{lemma}\label{lemma:ODE1}
	Let $\beta>0$ and $p\geq 0$ be real numbers, and let $\phibp(t)$
	be the solution of the Cauchy problem (\ref{defn:phibp-eq}),
	(\ref{defn:phibp-0}).
	
	Let $\ep$ and $K$ be positive constants, with $2\ep\beta\leq 1$,
	and let $G:[0,+\infty)\to[0,+\infty)$ be a function of class
	$C^{1}$ such that
	\begin{equation}
		G'(t)\leq -\frac{1}{\ep}\frac{1}{(1+t)^{p}}G(t)+
		\frac{K}{\ep}(1+t)^{p}\phibp(t)
		\quad\quad
		\forall t\geq 0.
		\label{hp:ODE1}
	\end{equation}
	
	Then we have that
	\begin{equation}
		G(t)\leq\left(2K+G(0)\right)(1+t)^{2p}\phibp(t)
		\quad\quad
		\forall t\geq 0.
		\label{th:ODE1}
	\end{equation}
\end{lemma}

\paragraph{\textit{\textmd{Proof}}}

Let us consider the differential equation
\begin{equation}
	y'(t)=-\frac{1}{\ep}\frac{1}{(1+t)^{p}}y(t)+
	\frac{K}{\ep}(1+t)^{p}\phibp(t)
	\quad\quad
	\forall t\geq 0.
	\label{eqn:ODE1}
\end{equation}

Assumption (\ref{hp:ODE1}) says that $G(t)$ is a subsolution of 
(\ref{eqn:ODE1}). Let $z(t)$ denote the right-hand side of 
(\ref{th:ODE1}). We claim that $z(t)$ is a supersolution of 
(\ref{eqn:ODE1}). Indeed a simple computation shows that
\begin{eqnarray*}
	z'(t) & = & 2p(2K+G(0))(1+t)^{2p-1}\phibp(t)+
	(2K+G(0))(1+t)^{2p}\phibp'(t)  \\
	\noalign{\vspace{1ex}}
	 & \geq & -\beta(2K+G(0))(1+t)^{p}\phibp(t)  \\
	 \noalign{\vspace{1ex}}
	 & \geq & -\frac{1}{\ep}(K+G(0))(1+t)^{p}\phibp(t)  \\
	 & = & -\frac{1}{\ep}\frac{1}{(1+t)^{p}}z(t)+
	\frac{K}{\ep}(1+t)^{p}\phibp(t),
\end{eqnarray*}
where in the second inequality we exploited that $2\ep\beta\leq 1$,
and $2G(0)\geq G(0)$.

Since $G(0)\leq z(0)$, estimate (\ref{th:ODE1}) follows from the 
standard comparison principle between subsolutions and 
supersolutions.\qed

\begin{lemma}\label{lemma:ODE2}
	Let $\psi_{1}:[0,+\infty)\to[0,+\infty)$ and 
	$\psi_{2}:[0,+\infty)\to[0,+\infty)$ be two continuous functions 
	such that
	$$K_{1}:=\int_{0}^{+\infty}\psi_{1}(t)\,dt<+\infty,
	\hspace{3em}
	K_{2}:=\int_{0}^{+\infty}\psi_{2}(t)\,dt<+\infty.$$
	
	Let $E:[0,+\infty)\to[0,+\infty)$ be a function of class $C^{1}$ 
	such that $E(0)=0$, and
	$$E'(t)\leq \psi_{1}(t)\sqrt{E(t)}+\psi_{2}(t)
		\quad\quad
		\forall t\geq 0.$$
	
	Then we have that
	\begin{equation}
		E(t)\leq K_{1}^{2}+2K_{2}
		\quad\quad
		\forall t\geq 0.
		\label{th:ODE2}
	\end{equation}
\end{lemma}

\paragraph{\textit{\textmd{Proof}}}

Let us fix any $T>0$. For every $t\in[0,T]$ we have that
$$E'(t)\leq \psi_{1}(t)\cdot
\biggl(\sup_{s\in[0,T]}E(s)\biggr)^{1/2}+
\psi_{2}(t).$$

Since $E(0)=0$, an easy integration gives that
$$E(t)\leq \biggl(\sup_{s\in[0,T]}E(s)\biggr)^{1/2}
\int_{0}^{t}\psi_{1}(s)\,ds
+\int_{0}^{t}\psi_{2}(s)\,ds\leq
K_{1}\biggl(\sup_{s\in[0,T]}E(s)\biggr)^{1/2}+K_{2}$$
for every $t\in[0,T]$.  Taking the supremum of the left-hand side as
$t\in[0,T]$, we obtain that
$$\sup_{s\in[0,T]}E(s) \leq K_{1}
\biggl(\sup_{s\in[0,T]}E(s)\biggr)^{1/2}+K_{2} \leq
\frac{1}{2}K_{1}^{2}+
\frac{1}{2}\biggl(\sup_{s\in[0,T]}E(s)\biggr)+K_{2},$$
hence 
$$\sup_{s\in[0,T]}E(s)\leq K_{1}^{2}+2K_{2},$$
and in particular $E(T)\leq K_{1}^{2}+2K_{2}$. Since $T$ is 
arbitrary, (\ref{th:ODE2}) is proved.\qed

\begin{lemma}\label{lemma:ODE3}
	Let $\beta>0$ and $p\geq 0$ be real numbers, and let $\phibp(t)$
	be the solution of the Cauchy problem (\ref{defn:phibp-eq}),
	(\ref{defn:phibp-0}).
	
	Let $\psi:[0,+\infty)\to[0,+\infty)$ be a continuous function 
	such that
	$$\int_{0}^{+\infty}\frac{\psi(s)}{\phibp(s)}\,ds<+\infty.$$
	
	Let $T>0$, and let $F:[T,+\infty)\to[0,+\infty)$ be a function of
	class $C^{1}$ such that
	\begin{equation}
		F'(t)\leq -\frac{\beta}{(1+t)^{p}}F(t)+\psi(t)
		\quad\quad
		\forall t\geq T.
		\label{hp:ODE3}
	\end{equation}
	
	Then we have that
	\begin{equation}
		F(t)\leq\left(\frac{F(T)}{\phibp(T)}+
		\int_{0}^{+\infty}\frac{\psi(s)}{\phibp(s)}\,ds\right)
		\cdot\phibp(t)
		\quad\quad
		\forall t\geq T.
		\label{th:ODE3}
	\end{equation}
\end{lemma}

\paragraph{\textit{\textmd{Proof}}}

Let us consider the differential equation
\begin{equation}
	y'(t)=-\frac{\beta}{(1+t)^{p}}y(t)+\psi(t)
	\quad\quad
	\forall t\geq 0.
	\label{eqn:ODE3}
\end{equation}

Assumption (\ref{hp:ODE3}) says that $F(t)$ is a subsolution of 
(\ref{eqn:ODE3}) for $t\geq T$. On the other hand, it is easy to see 
that
$$z(t):=\left(\frac{F(T)}{\phibp(T)}+
\int_{T}^{t}\frac{\psi(s)}{\phibp(s)}\,ds\right)
\cdot\phibp(t)$$
is a solution of (\ref{eqn:ODE3}) for $t\geq T$.  Since $F(T)=z(T)$,
the standard comparison principle between subsolutions and
supersolutions implies that $F(t)\leq z(t)$ for every $t\geq T$, which
in turn implies (\ref{th:ODE3}).\qed

\subsection{Proof of Theorem~\ref{thm:lin-h}}

Let us describe the strategy of the proof before entering into
details.  Let us take any admissible value $\beta$, which means any
$\beta\in(0,2\mu\nu)$ if $p=0$, and any $\beta>0$ if $p>0$.  Let
$\phibp(t)$ be the solution of the Cauchy problem
(\ref{defn:phibp-eq}), (\ref{defn:phibp-0}).
	
Estimates (\ref{th:h-p=0}) through (\ref{th:h-p=1}) are equivalent to 
showing that
\begin{equation}
	\Gamma_{\ep}(t)\leq k_{1}\phibp(t)
	\quad\quad
	\forall t\geq 0
	\label{th:phibp}
\end{equation}
for the admissible values of $\beta$.

Let $\mu$ be the constant in (\ref{hp:cep-below}), and let us choose 
$\delta$ and $T$ in such a way that
\begin{equation}
	\delta:=\frac{2(\beta+1)\nu}{2\mu\nu-\beta},
	\hspace{3em}
	T:=0
	\label{choice:h-p=0}
\end{equation}
if $p=0$ (note that $\delta>0$), and
\begin{equation}
	\delta:=\frac{\beta+2}{\mu},
	\hspace{3em}
	(1+T)^{2p}\geq\frac{\delta\beta}{2\nu}
	\label{choice:h-p>0}
\end{equation}
if $p>0$.
For every $\ep\in(0,\ep_{0})$, we consider the energies
\begin{equation}
	\Eep(t):=\frac{\ep|\uep'(t)|^{2}}{\cep(t)}+|A^{1/2}\uep(t)|^{2},
	\label{defn:Eep}
\end{equation}
\begin{equation}
	\Fep(t):=\frac{\ep|\uep'(t)|^{2}}{\cep(t)}+|A^{1/2}\uep(t)|^{2}
	+\frac{\ep\delta}{(1+t)^{p}}\langle\uep'(t),\uep(t)\rangle
	+\frac{\delta}{2}\frac{1}{(1+t)^{2p}}|\uep(t)|^{2}.
	\label{defn:Fep}
\end{equation}

We claim that there exist $\ep_{2}\in(0,\ep_{0})$, and positive
constants $k_{2}$, \ldots, $k_{5}$, such that
\begin{equation}
	k_{2}\left(\ep|\uep'(t)|^{2}+|A^{1/2}\uep(t)|^{2}\right)\leq
	\Eep(t)\leq 
	k_{3}\left(\ep|\uep'(t)|^{2}+|A^{1/2}\uep(t)|^{2}\right),
	\label{est:Eep-equiv}
\end{equation}
\begin{equation}
	k_{4}\left(\ep|\uep'(t)|^{2}+|A^{1/2}\uep(t)|^{2}\right)\leq
	\Fep(t)\leq 
	k_{5}\left(\ep|\uep'(t)|^{2}+|A^{1/2}\uep(t)|^{2}\right)
	\label{est:Fep-equiv}
\end{equation}
for every $t\geq 0$ and every $\ep\in(0,\ep_{2})$. 
Moreover we claim that
\begin{equation}
	\Eep'(t)\leq 0
	\quad\quad
	\forall t\geq 0,\ \forall\ep\in(0,\ep_{2}),
	\label{est:Eep'}
\end{equation}
\begin{equation}
	\Fep'(t)\leq -\frac{\beta}{(1+t)^{p}}\Fep(t)
	\quad\quad
	\forall t\geq T,\ \forall\ep\in(0,\ep_{2}).
	\label{est:Fep'}
\end{equation}

Let us assume that we have proved these claims.  Thanks to
(\ref{est:Eep'}), and to the estimate from below in
(\ref{est:Eep-equiv}), we have that
$$\ep|\uep'(t)|^{2}+|A^{1/2}\uep(t)|^{2}\leq
\frac{1}{k_{2}}\Eep(t)\leq\frac{1}{k_{2}}\Eep(0)\leq k_{6}$$
for every $t\geq 0$. Since $\phibp(t)$ is decreasing, this implies that
\begin{equation}
	\ep|\uep'(t)|^{2}+|A^{1/2}\uep(t)|^{2}\leq
	\frac{k_{6}}{\phibp(T)}\cdot\phibp(t)= k_{7}\phibp(t)
	\quad\quad
	\forall t\in[0,T].
	\label{est:0-T}
\end{equation}

For $t\geq T$, we exploit (\ref{est:Fep'}).  First of all, from
(\ref{est:0-T}) with $t=T$, and the estimate from above in
(\ref{est:Fep-equiv}), we have that
$$\Fep(T)\leq
k_{5}\left(\ep|\uep'(T)|^{2}+|A^{1/2}\uep(T)|^{2}\right)\leq
k_{8}\phibp(T).$$

Therefore, from Lemma~\ref{lemma:ODE3} applied with $\psi(t)\equiv 0$,
we deduce that $\Fep(t)\leq k_{8}\phibp(t)$ for every $t\geq T$.
Exploiting this inequality, the estimate from below in
(\ref{est:Fep-equiv}), and (\ref{est:0-T}), we conclude that
\begin{equation}
	\ep|\uep'(t)|^{2}+|A^{1/2}\uep(t)|^{2}\leq k_{9}\phibp(t)
	\quad\quad
	\forall t\geq 0,\ \forall\ep\in(0,\ep_{2}).
	\label{est:order-0}
\end{equation}

Since the operator is coercive, this estimate on
$|A^{1/2}\uep(t)|^{2}$ yields an analogous estimate on
$|\uep(t)|^{2}$.

Up to now, we only assumed that $(u_{0},u_{1})\in D(A^{1/2})\times H$.
Let us assume now that $(u_{0},u_{1})\in D(A)\times D(A^{1/2})$.
Since equation (\ref{pbm:h-leq}) is linear, estimate
(\ref{est:order-0}) can be applied to $A^{1/2}\uep(t)$, which is once
again a solution to (\ref{pbm:h-leq}).  We thus obtain that
\begin{equation}
	\ep|A^{1/2}\uep'(t)|^{2}+|A\uep(t)|^{2}\leq k_{10}\phibp(t)
	\quad\quad
	\forall t\geq 0,\ \forall\ep\in(0,\ep_{2}).
	\label{est:order-1}
\end{equation}

It remains to prove the $\ep$-independent estimate on
$|\uep'(t)|^{2}$.  To this end, we set
\begin{equation}
	\Gep(t):=|\uep'(t)|^{2},
	\label{defn:Gep}
\end{equation}
and we claim that
\begin{equation}
	\Gep'(t)\leq -\frac{1}{\ep}\frac{1}{(1+t)^{p}}\Gep(t)+
	\frac{k_{11}}{\ep}(1+t)^{p}\phibp(t)
	\quad\quad
	\forall t\geq 0,\ \forall\ep\in(0,\ep_{2}).
	\label{est:Gep'}
\end{equation}

If we prove the claim, then from Lemma~\ref{lemma:ODE1} if follows that
\begin{equation}
	|\uep'(t)|^{2}=\Gep(t)\leq
	k_{12}(1+t)^{2p}\phibp(t)
	\quad\quad
	\forall t\geq 0,\ \forall\ep\in(0,\ep_{2}).
	\label{est:uep'-p}
\end{equation}

What we actually need is the same estimate without the factor
$(1+t)^{2p}$.  If $p=0$, there is nothing to do.  If $p>0$, we take
$\beta'=\beta+2$, and from (\ref{est:uep'-p}) we obtain that
$$|\uep'(t)|^{2}=\Gep(t)\leq
k_{13}(1+t)^{2p}\Phi_{\beta',p}(t) \quad\quad
\forall t\geq 0,\ \forall\ep\in(0,\ep_{1}),$$
of course with new positive constants $k_{13}$ and $\ep_{1}\leq\ep_{2}$,
depending also on $\beta'$.

Finally, our choice of $\beta'$ guarantees that
$$(1+t)^{2p}\Phi_{\beta',p}(t)\leq k_{14}\phibp(t) \quad\quad
\forall t\geq 0$$
for a suitable $k_{14}$ depending on $p$, $\beta$, $\beta'$ (this
inequality can be easily proved exploiting the explicit formulae for
$\phibp(t)$ and $\Phi_{\beta',p}(t)$, and the fact that $p\leq
1$).  This completes the proof of (\ref{th:phibp}) for every 
$\ep\in(0,\ep_{1})$.

So we are left to proving (\ref{est:Eep-equiv}) through 
(\ref{est:Fep'}), and (\ref{est:Gep'}).

\paragraph{\textmd{\textit{Equivalence between energies}}}

Due to (\ref{hp:cep-below}) and (\ref{hp:cep-above}), estimate
(\ref{est:Eep-equiv}) holds true with
$$k_{2}:=\min\left\{\frac{1}{M_{3}},1\right\},
\hspace{3em}
k_{3}:=\max\left\{\frac{1}{\mu},1\right\}.$$

In order to prove (\ref{est:Fep-equiv}), let us estimate separately 
the four terms in (\ref{defn:Fep}). Due to (\ref{hp:cep-below}) and 
(\ref{hp:cep-above}), we have that
$$\frac{\ep|\uep'(t)|^{2}}{M_{3}}\leq
\frac{\ep|\uep'(t)|^{2}}{\cep(t)}\leq
\frac{\ep|\uep'(t)|^{2}}{\mu}.$$

Due to (\ref{hp:coercive}) we have that
$$0\leq\frac{\delta}{2}\frac{1}{(1+t)^{2p}}|\uep(t)|^{2}\leq
\frac{\delta}{2}|\uep(t)|^{2}\leq
\frac{\delta}{2\nu}|A^{1/2}\uep(t)|^{2}.$$

Applying once again (\ref{hp:coercive}), and the inequality between
arithmetic and geometric mean, we obtain that
$$\frac{\ep|\uep'(t)|^{2}}{2M_{3}}+\frac{1}{2}|A^{1/2}\uep(t)|^{2}
\geq \frac{\ep|\uep'(t)|^{2}}{2M_{3}}+\frac{\nu}{2}|\uep(t)|^{2}
\geq\sqrt{\frac{\ep\nu}{M_{3}}}\cdot|\uep'(t)|\cdot|\uep(t)|.$$

If $\ep$ is small enough, this implies that
$$\frac{\ep|\uep'(t)|^{2}}{2M_{3}}+\frac{1}{2}|A^{1/2}\uep(t)|^{2}\geq
\frac{\ep\delta}{(1+t)^{2p}}
\left|\langle\uep'(t),\uep(t)\rangle\right|.$$

From all these estimates, we easily obtain that 
$$\Fep(t)\geq \frac{\ep|\uep'(t)|^{2}}{2M_{3}}+
\frac{1}{2}|A^{1/2}\uep(t)|^{2},$$
and 
$$\Fep(t)\leq \frac{\ep|\uep'(t)|^{2}}{\mu}+
|A^{1/2}\uep(t)|^{2}+
\frac{\delta}{2\nu}|A^{1/2}\uep(t)|^{2}+
\frac{\ep|\uep'(t)|^{2}}{2M_{3}}+
\frac{1}{2}|A^{1/2}\uep(t)|^{2},$$
from which (\ref{est:Fep-equiv}) follows with 
$$k_{4}:=\min\left\{\frac{1}{2M_{3}},\frac{1}{2}\right\},
\hspace{3em}
k_{5}:=\max\left\{\frac{1}{\mu}+\frac{1}{2M_{3}},
\frac{3}{2}+\frac{\delta}{2\nu}\right\}.$$

\paragraph{\textmd{\textit{Differential inequality for $\Eep$}}}

The time-derivative of (\ref{defn:Eep}) is
$$\Eep'(t)=-\frac{1}{(1+t)^{p}}\frac{|\uep'(t)|^{2}}{\cep(t)}
\left(2+\ep\frac{\cep'(t)(1+t)^{p}}{\cep(t)}\right).$$

From (\ref{hp:cep-below}) and (\ref{hp:cep'-above}) we have that
$$\ep\frac{|\cep'(t)|(1+t)^{p}}{\cep(t)}\leq
\frac{M_{4}}{\mu}\ep,$$
so that $\Eep'(t)\leq 0$ for every $t\geq 0$, provided that $\ep$ is
small enough.  This proves (\ref{est:Eep'}).

\paragraph{\textmd{\textit{Differential inequality for $\Fep$}}}

The time-derivative of (\ref{defn:Fep}) is
\begin{eqnarray*}
	\Fep'(t) & = & -\frac{1}{(1+t)^{p}}\frac{|\uep'(t)|^{2}}{\cep(t)}
	\left(2+\ep\frac{\cep'(t)(1+t)^{p}}{\cep(t)}-
	\ep\delta\cep(t)\right)   \\
	 &  & - \frac{\delta\cep(t)}{(1+t)^{p}}|A^{1/2}\uep(t)|^{2} 
	 -\delta p\frac{|\uep(t)|^{2}}{(1+t)^{2p+1}}
	 -\frac{\ep\delta p}{(1+t)^{1+p}}\langle\uep'(t),\uep(t)\rangle
\end{eqnarray*}

Therefore (\ref{est:Fep'}) holds true if and only if
$$\hspace{-3em}
\dfrac{|\uep'(t)|^{2}}{\cep(t)}
\left(2+\ep\dfrac{\cep'(t)(1+t)^{p}}{\cep(t)}-
\ep\delta\cep(t)-\ep\beta\right)+
(\delta\cep(t)-\beta)|A^{1/2}\uep(t)|^{2}+$$
\begin{equation}
	\left(\dfrac{\delta p}{(1+t)^{1+p}}
	-\dfrac{\delta\beta}{2}\dfrac{1}{(1+t)^{2p}}\right)|\uep(t)|^{2}
	+\left(\dfrac{\ep\delta p}{1+t}
	-\dfrac{\ep\delta\beta}{(1+t)^{p}}\right)
	\langle\uep'(t),\uep(t)\rangle\geq 0
	\label{est:mostro}
\end{equation}
holds true for every $t\geq T$, and every $\ep$ small enough.

Let $S_{1}$, \ldots, $S_{4}$ denote the four terms in 
(\ref{est:mostro}). Due to (\ref{hp:cep-below}) through 
(\ref{hp:cep'-above}), for every small enough $\ep$ we have that
$$\ep\dfrac{|\cep'(t)|(1+t)^{p}}{\cep(t)}\leq
\frac{M_{4}}{\mu}\ep\leq\frac{1}{3},
\quad\quad
\ep\delta\cep(t)\leq\ep\delta M_{3}\leq\frac{1}{3},
\quad\quad
\ep\beta\leq\frac{1}{3},$$
hence
\begin{equation}
	S_{1}\geq\frac{|\uep'(t)|^{2}}{\cep(t)}\geq
	\frac{1}{M_{3}}|\uep'(t)|^{2}.
	\label{est:S1}
\end{equation}

Since $\delta\mu\geq\beta$, from (\ref{hp:coercive}) we have that
\begin{eqnarray*}
	S_{2}+S_{3} & \geq & (\delta \mu-\beta)|A^{1/2}\uep(t)|^{2}-
	\frac{\delta\beta}{2}\frac{1}{(1+t)^{2p}}|\uep(t)|^{2} \\
	 & \geq & \left[(\delta \mu-\beta)\nu-\frac{\delta\beta}{2}
	\frac{1}{(1+T)^{2p}}\right]|\uep(t)|^{2}
\end{eqnarray*}
for every $t\geq T$.  Due to the choices (\ref{choice:h-p=0}) and
(\ref{choice:h-p>0}), in both cases the term in brackets is greater
than or equal to $\nu$, hence $S_{2}+S_{3}\geq\nu |\uep(t)|^{2}$ for
every $t\geq T$.  Adding this inequality to (\ref{est:S1}), and
applying the inequality between arithmetic and geometric mean, we
deduce that $$S_{1}+S_{2}+S_{3} \geq
\frac{1}{M_{3}}|\uep'(t)|^{2}+\nu|\uep(t)|^{2} \geq
2\sqrt{\frac{\nu}{M_{3}}}\cdot|\uep'(t)|\cdot|\uep(t)|.$$

As a consequence, if $\ep$ is small enough and $t\geq T$, we have that
$$S_{1}+S_{2}+S_{3} \geq\ep\delta(1+\beta)|\uep'(t)|\cdot|\uep(t)|
\geq \left(\frac{\ep\delta p}{1+t}
+\frac{\ep\delta\beta}{(1+t)^{p}}\right) |\uep'(t)|\cdot|\uep(t)| \geq
|S_{4}|,$$
which proves (\ref{est:mostro}), hence also (\ref{est:Fep'}).

\paragraph{\textmd{\textit{Differential inequality for $\Gep$}}}

The time-derivative of (\ref{defn:Gep}) is
$$\Gep'(t)=-\frac{2}{\ep}\frac{1}{(1+t)^{p}}|\uep'(t)|^{2}
-\frac{2}{\ep}\cep(t)\langle A\uep(t),\uep'(t)\rangle.$$

From (\ref{hp:cep-above}) we have that
$$-2\cep(t)\langle A\uep(t),\uep'(t)\rangle\leq
2M_{3}|\uep'(t)|\cdot|A\uep(t)|\leq
\frac{|\uep'(t)|^{2}}{(1+t)^{p}}+M_{3}^{2}(1+t)^{p}|A\uep(t)|^{2},$$
hence
$$\Gep'(t)\leq -\frac{1}{\ep}\frac{1}{(1+t)^{p}}\Gep(t)+
\frac{M_{3}^{2}}{\ep}(1+t)^{p}|A\uep(t)|^{2}.$$

At this point (\ref{est:Gep'}) follows from (\ref{est:order-1}).  

The proof of Theorem~\ref{thm:lin-h} is thus complete.\qed

\subsection{Singular perturbation: preliminary estimates}

In this subsection we begin the analysis of the singular perturbation 
problem in the linear setting. If we set
\begin{equation}
	\gep(t):=-(\cep(t)-c(t))Au(t)-\ep u''(t),
	\label{defn:gep}
\end{equation}
we have that $\rep(t)$ and $\roep(t)$ satisfy
\begin{equation}
	\ep\rep''(t)+\frac{1}{(1+t)^{p}}\rep'(t)+\cep(t)A\roep(t)=\gep(t),
	\label{eqn:rep-lin}
\end{equation}
and 
$$\roep(0)=0,
	\quad\quad
	\rep'(0)=0.$$

In the next two results we prove estimates on $\gep(t)$ and on the 
corrector $\tetep(t)$. 

\begin{lemma}
	Let us consider the same assumptions of Theorem~\ref{thm:lin-sp}.
	Let $\gep(t)$ be defined according to (\ref{defn:gep}).  Let
	$\phibp(t)$ be the solution of the Cauchy problem
	(\ref{defn:phibp-eq}), (\ref{defn:phibp-0}), with $\beta>0$ if
	$p>0$, and $0<\beta<2\mu\nu$ if $p=0$.
	
	Then we have the following estimates.
	\begin{enumerate}
		\renewcommand{\labelenumi}{(\arabic{enumi})} 
		\item If $u_{0}\in D(A^{3/2})$, then there exists a constant
		$C$ such that
		\begin{equation}
			\int_{0}^{+\infty}\frac{(1+t)^{p}}{\phibp(t)}
			\cdot|\gep(t)|^{2}\,dt
			\leq C\ep^{2}
			\quad\quad
			\forall\ep\in(0,\ep_{0}).
			\label{th:gep-int-1}
		\end{equation}
	
		\item If in addition we have that $u_{0}\in D(A^{2})$, then
		there exists a constant $C$ such that
		\begin{equation}
			\int_{0}^{+\infty}\frac{(1+t)^{p}}{\phibp(t)}
			\cdot|A^{1/2}\gep(t)|^{2}\,dt
			\leq C\ep^{2}
			\quad\quad
			\forall\ep\in(0,\ep_{0}),
			\label{th:gep-int-2}
		\end{equation}
		\begin{equation}
			|\gep(t)|^{2}\leq C\ep^{2}\phibp(t)
			\quad\quad
			\forall t\geq 0,\ \forall\ep\in(0,\ep_{0}).
			\label{th:gep-pointwise}
		\end{equation}

	\end{enumerate}
\end{lemma}

\paragraph{\textit{\textmd{Proof}}}

From (\ref{defn:gep}) and (\ref{hp:cep-c}) we have that
$$|\gep(t)|^{2}\leq k_{1}\ep^{2}|Au(t)|^{2}+2\ep^{2}|u''(t)|^{2}.$$

We can estimate $|Au(t)|^{2}$ and $|u''(t)|^{2}$, or their integrals,
by means of Proposition~\ref{prop:p-lin}.  To this end, let us
consider the function $\psiap(t)$ defined in (\ref{defn:psiap}).  We
claim that, for every admissible value of $p$ and $\beta$, there
exists $\alpha>0$ for which Proposition~\ref{prop:p-lin} applies, and
such that
\begin{equation}
	\frac{(1+t)^{p}}{\phibp(t)}\leq\frac{k_{2}}{\psiap(t)}
	\quad\quad
	\forall t\geq 0.
	\label{est:phibp-psiap}
\end{equation}

Indeed it is enough to take $\alpha=\beta$ if $p=0$ (in which case
there is basically nothing to prove), and any $\alpha\in(0,\gamma)$ if
$p>0$ (because in this case $\psiap(t)$ has an exponential decay rate
which is faster than the decay rate of $\phibp(t)$).  Thus we have
that
$$\int_{0}^{+\infty}\frac{(1+t)^{p}}{\phibp(t)} \cdot|\gep(t)|^{2}\,dt
\leq k_{3}\ep^{2}\left(\int_{0}^{+\infty}
\frac{|Au(t)|^{2}}{\psiap(t)}\,dt+\int_{0}^{+\infty}
\frac{|u''(t)|^{2}}{\psiap(t)}\,dt\right),$$
so that (\ref{th:gep-int-1}) follows from (\ref{th:p-ak-int}) with 
$k=1$, and (\ref{th:u''-int-1}).

The proof of (\ref{th:gep-int-2}) is analogous: we just exploit
(\ref{th:p-ak-int}) with $k=2$, and (\ref{th:u''-int-2}) instead of
(\ref{th:u''-int-1}).

It remains to prove (\ref{th:gep-pointwise}). Let $\gamma$ be the 
constant defined in (\ref{defn:psiap}). Then, in analogy with 
(\ref{est:phibp-psiap}), we have that
$$\frac{(1+t)^{4p}}{\phibp(t)}\leq\frac{k_{4}}{\psigp(t)}
\quad\quad
\forall t\geq 0,$$
hence
$$\frac{|\gep(t)|^{2}}{\phibp(t)}=\frac{(1+t)^{4p}}{\phibp(t)}\cdot
|\gep(t)|^{2}\cdot\frac{1}{(1+t)^{4p}}\leq
k_{4}\frac{|\gep(t)|^{2}}{\psigp(t)}\cdot\frac{1}{(1+t)^{4p}}$$
$$\leq k_{5}\ep^{2}\frac{|Au(t)|^{2}}{\psigp(t)}+
k_{6}\ep^{2}\frac{|u''(t)|^{2}}{\psigp(t)}
\cdot\frac{1}{(1+t)^{4p}}.$$

At this point (\ref{th:gep-pointwise}) follows from (\ref{th:p-ak}) with 
$k=2$, and (\ref{th:u''-pointwise}).\qed

\begin{lemma}
	Let us consider the same assumptions of Theorem~\ref{thm:lin-sp}.
	Let $\tetep(t)$ be the solution of the Cauchy problem
	(\ref{pbm:tetep-eq}), (\ref{pbm:tetep-ldata}).  Let $\phibp(t)$ be
	the solution of the Cauchy problem (\ref{defn:phibp-eq}),
	(\ref{defn:phibp-0}).
	
	Let us assume that $4\ep_{0}\leq 1$, $2\ep_{0}\beta\leq 1$, and 
	that $\theta_{0}\in D(A^{(k+1)/2})$ for some $k\in\n$.
	
	Then there exists a constant $C$ such that for every 
	$\ep\in(0,\ep_{0})$ we have that
	\begin{equation}
		\int_{0}^{+\infty}\frac{1}{\phibp(t)}\cdot\left(
		|A^{k/2}\tetep'(t)|+|A^{k/2}\tetep'(t)|^{2}+
		|A^{(k+1)/2}\tetep'(t)|\right)\,dt\leq C\ep.
		\label{th:tetep'-int}
	\end{equation}
	
\end{lemma}

\paragraph{\textit{\textmd{Proof}}}

Let $\zep(t)$ be the solution of equation
\begin{equation}
	\ep\zep'(t)+\frac{1}{(1+t)^{p}}\zep(t)=0
	\quad\quad
	\forall t\geq 0,
	\label{defn:zep}
\end{equation}
with initial condition $\zep(0)=1$. It is easy to see that 
$\tetep'(t)=\theta_{0}\zep(t)$.

Since $0\leq\zep(t)\leq 1$ for every $t\geq 0$, we have also that 
$\zep^{2}(t)\leq\zep(t)$. Therefore, (\ref{th:tetep'-int}) is proved if we 
show that
\begin{equation}
	\int_{0}^{+\infty}\frac{\zep(t)}{\phibp(t)}\,dt\leq4\ep.
	\label{est:zep}
\end{equation}

Let us set $\wep(t):=\zep(t)\cdot[\phibp(t)]^{-1}$. From (\ref{defn:zep}) 
and (\ref{defn:phibp-eq}), it turns out that $\wep(t)$ is the 
solution of the ordinary differential equation
\begin{equation}
	\wep'(t)=-\left(\frac{1}{\ep}-\beta\right)
	\frac{1}{(1+t)^{p}}\wep(t)
	\quad\quad
	\forall t\geq 0,
	\label{defn:wep}
\end{equation}
with initial datum $\wep(0)=1$. On the other hand, when 
$2\ep\beta\leq 1$, it is easy to show that 
$\yep(t):=(1+t)^{-1/(2\ep)}$ is a supersolution of (\ref{defn:wep}). 
Indeed we have that
$$\yep'(t)=-\frac{1}{2\ep}\frac{\yep(t)}{1+t}\geq
-\frac{1}{2\ep}\frac{\yep(t)}{(1+t)^{p}}\geq
-\left(\frac{1}{\ep}-\beta\right)\frac{\yep(t)}{(1+t)^{p}}.$$

Since $\yep(0)=\wep(0)$, the standard comparison principle gives that 
$\wep(t)\leq\yep(t)$ for every $t\geq 0$. Since $4\ep\leq 1$, it 
follows that
$$\int_{0}^{+\infty}\wep(t)\,dt\leq
\int_{0}^{+\infty}\frac{1}{(1+t)^{1/(2\ep)}}\,dt=
\frac{2\ep}{1-2\ep}\leq 4\ep.$$

This completes the proof of (\ref{est:zep}), hence also the proof of 
(\ref{th:tetep'-int}).\qed

\subsection{Proof of Theorem~\ref{thm:lin-sp}}

Let us describe the strategy of the proof, which is similar to
Theorem~\ref{thm:lin-h}.  Let us take any admissible value $\beta$,
which means any $\beta\in(0,2\mu\nu)$ if $p=0$, and any $\beta>0$ if
$p>0$.  Let $\phibp(t)$ be the solution of the Cauchy problem
(\ref{defn:phibp-eq}), (\ref{defn:phibp-0}).

The conclusions of statement~(1) of Theorem~\ref{thm:lin-sp} are
equivalent to showing that
\begin{equation}
	\Gamma_{r,\ep}(t)\leq k_{1}\phibp(t)
	\quad\quad
	\forall t\geq 0
	\label{th:sp-phibp}
\end{equation}
for the admissible values of $\beta$.

Let $\mu$ be the constant in (\ref{hp:cep-below}), and let us choose
$\delta$, $\sigma$, $T$ in such a way that
\begin{equation}
	\delta:=\frac{4(\beta+1)\nu}{2\mu\nu-\beta},
	\hspace{3em}
	\sigma:=\mu\nu-\frac{\beta}{2},
	\hspace{3em}
	T:=0
	\label{choice:sp-p=0}
\end{equation}
if $p=0$ (note that $\delta>0$), and
\begin{equation}
	\delta:=\frac{\beta+2}{\mu},
	\hspace{3em}
	\sigma:=1,
	\hspace{3em}
	(1+T)^{2p}\geq\frac{\delta}{2\nu}(\beta+\sigma)
	\label{choice:sp-p>0}
\end{equation}
if $p>0$.

For every $\ep\in(0,\ep_{0})$, we consider the energies
\begin{equation}
	\EEep(t):=\frac{\ep|\rep'(t)|^{2}}{\cep(t)}+|A^{1/2}\roep(t)|^{2},
	\label{defn:EEep}
\end{equation}
\begin{equation}
	\FFep(t):=\frac{\ep|\rep'(t)|^{2}}{\cep(t)}+|A^{1/2}\roep(t)|^{2}
	+\frac{\ep\delta}{(1+t)^{p}}\langle\rep'(t),\roep(t)\rangle
	+\frac{\delta}{2}\frac{1}{(1+t)^{2p}}|\roep(t)|^{2}.
	\label{defn:FFep}
\end{equation}

The arguments used in the proof of (\ref{est:Eep-equiv}) and
(\ref{est:Fep-equiv}) can be adapted word-by-word to the energies
$\EEep(t)$ and $\FFep(t)$.  We obtain that there exist positive
constants $k_{2}$, \ldots, $k_{5}$ such that
\begin{equation}
	k_{2}\left(\ep|\rep'(t)|^{2}+|A^{1/2}\roep(t)|^{2}\right)\leq
	\EEep(t)\leq 
	k_{3}\left(\ep|\rep'(t)|^{2}+|A^{1/2}\roep(t)|^{2}\right),
	\label{est:EEep-equiv}
\end{equation}
\begin{equation}
	k_{4}\left(\ep|\rep'(t)|^{2}+|A^{1/2}\roep(t)|^{2}\right)\leq
	\FFep(t)\leq 
	k_{5}\left(\ep|\rep'(t)|^{2}+|A^{1/2}\roep(t)|^{2}\right)
	\label{est:FFep-equiv}
\end{equation}
for every $t\geq 0$, provided that $\ep$ is small enough.

Moreover, we claim that there exists $\ep_{2}\in(0,\ep_{0})$ such that
\begin{equation}
	\EEep'(t)\leq \psi_{1,\ep}(t)\sqrt{\EEep(t)}+\psi_{2,\ep}(t)
	\quad\quad
	\forall t\geq 0,\ \forall\ep\in(0,\ep_{2}),
	\label{est:EEep'}
\end{equation}
\begin{equation}
	\FFep'(t)\leq -\frac{\beta}{(1+t)^{p}}\FFep(t)+\psi_{3,\ep}(t)
	\quad\quad
	\forall t\geq T,\ \forall\ep\in(0,\ep_{2}),
	\label{est:FFep'}
\end{equation}
where the functions $\psi_{i,\ep}(t)$ (with $i=1,2,3$) are nonnegative
continuous functions depending on $A^{1/2}\tetep'(t)$ and $\gep(t)$,
and such that
\begin{equation}
	\int_{0}^{+\infty}\psi_{1,\ep}(t)\,dt\leq k_{6}\ep,
	\hspace{3em}
	\int_{0}^{+\infty}\psi_{2,\ep}(t)\,dt\leq k_{7}\ep^{2},
	\label{hp:psi12}
\end{equation}
\begin{equation}
	\int_{0}^{+\infty}\frac{\psi_{3,\ep}(t)}{\phibp(t)}\,dt\leq k_{8}\ep^{2}.
	\label{hp:psi3}
\end{equation}

Let us assume that we have proved these claims.  Thanks to 
(\ref{est:EEep'}) and (\ref{hp:psi12}), we can apply 
Lemma~\ref{lemma:ODE2} to the function $\EEep(t)$ (note that now 
$\EEep(0)=0$). We obtain that
\begin{equation}
	\EEep(t)\leq k_{9}\ep^{2}
	\quad\quad
	\forall t\geq 0.
	\label{est:sp-EEep}
\end{equation}

Due to the estimate from below in
(\ref{est:EEep-equiv}), this implies that 
$$\ep|\rep'(t)|^{2}+|A^{1/2}\roep(t)|^{2}\leq
\frac{1}{k_{2}}\EEep(t)\leq k_{10}\ep^{2}$$
for every $t\geq 0$.  Since $\phibp(t)$ is decreasing, we can conclude
that
\begin{equation}
	\ep|\rep'(t)|^{2}+|A^{1/2}\roep(t)|^{2}\leq 
	\frac{k_{10}\ep^{2}}{\phibp(T)}\cdot\phibp(t)= k_{11}\ep^{2}\phibp(t)
	\quad\quad
	\forall t\in[0,T].
	\label{est:sp-0-T}
\end{equation}

For $t\geq T$, we exploit (\ref{est:FFep'}).  First of all, from
(\ref{est:sp-0-T}) with $t=T$, and the estimate from above in
(\ref{est:FFep-equiv}), we have that
$$\FFep(T)\leq
k_{5}\left(\ep|\rep'(T)|^{2}+|A^{1/2}\roep(T)|^{2}\right)\leq
k_{12}\ep^{2}\phibp(T).$$

Due to (\ref{est:FFep'}) and (\ref{hp:psi3}), we can apply
Lemma~\ref{lemma:ODE3} to the function $\FFep(t)$.  We obtain that
$\FFep(t)\leq k_{13}\ep^{2}\phibp(t)$ for every $t\geq T$.  Exploiting
this inequality, the estimate from below in (\ref{est:FFep-equiv}),
and (\ref{est:sp-0-T}), we conclude that
$$\ep|\rep'(t)|^{2}+|A^{1/2}\roep(t)|^{2}\leq k_{14}\ep^{2}\phibp(t)
	\quad\quad
	\forall t\geq 0,\ \forall\ep\in(0,\ep_{2}).$$

Since the operator is coercive, this estimate on
$|A^{1/2}\roep(t)|^{2}$ yields an analogous estimate on
$|\roep(t)|^{2}$.  This completes the proof of (\ref{th:sp-phibp}),
hence of statement~(1), for initial data $(u_{0},u_{1})\in
D(A^{3/2})\times D(A^{1/2})$, the regularity of data being required in
the verification of (\ref{hp:psi12}) and (\ref{hp:psi3}).

Let us proceed now to statement~(2), where it is assumed that 
$(u_{0},u_{1})\in D(A^{2})\times D(A)$, and it is required to prove 
in addition that
$$\ep|A^{1/2}\rep'(t)|^{2}+|A\roep(t)|^{2}+|\rep'(t)|^{2}
	\leq k_{15}\ep^{2}\phibp(t)
	\quad\quad
	\forall t\geq 0,\ \forall\ep\in(0,\ep_{1})$$
for some $\ep_{1}\in(0,\ep_{2}]$.  Due to the linearity of
(\ref{eqn:rep-lin}), an analogous identity holds true with
$A^{1/2}\roep(t)$, $A^{1/2}\rep(t)$, and $A^{1/2}\gep(t)$ instead of
$\roep(t)$, $\rep(t)$, and $\gep(t)$, respectively.  So we can repeat
the arguments used sofar, paying attention to verifying
(\ref{hp:psi12}) and (\ref{hp:psi3}) also for the new functions
$\psi_{\ep,i}(t)$, which now depend on $A\tetep'(t)$ and
$A^{1/2}\gep(t)$.  We end up with
\begin{equation}
	\ep|A^{1/2}\rep'(t)|^{2}+|A\roep(t)|^{2}\leq k_{16}\ep^{2}\phibp(t)
	\quad\quad
	\forall t\geq 0,\ \forall\ep\in(0,\ep_{2}).
	\label{est:sp-order-1}
\end{equation}

It remains to prove the $\ep$-independent estimates on $\rep'(t)$. To 
this end, we set
\begin{equation}
	\GGep(t):=|\rep'(t)|^{2},
	\label{defn:GGep}
\end{equation}
and we claim that
\begin{equation}
	\GGep'(t)\leq -\frac{1}{\ep}\frac{1}{(1+t)^{p}}\GGep(t)+
	\frac{1}{\ep}(1+t)^{p}\cdot k_{17}\ep^{2}\phibp(t)
	\quad\quad
	\forall t\geq 0,\ \forall\ep\in(0,\ep_{2}).
	\label{est:GGep'}
\end{equation}

If we prove the claim, then from Lemma~\ref{lemma:ODE1} it follows 
that (note that now $\GGep(0)=0$)
$$|\rep'(t)|^{2}=\GGep(t)\leq
	k_{18}\ep^{2}(1+t)^{2p}\phibp(t)
	\quad\quad
	\forall t\geq 0,\ \forall\ep\in(0,\ep_{2}).$$

Finally, when $p>0$, we can get free of the factor $(1+t)^{2p}$
exactly as in the proof of Theorem~\ref{thm:lin-h}, possibly changing
$\ep_{2}$ with some smaller $\ep_{1}$.

So we are left to proving (\ref{est:EEep'}) through (\ref{hp:psi3}),
both in the case of initial data $(u_{0},u_{1})\in D(A^{3/2})\times
D(A^{1/2})$, and in the case $(u_{0},u_{1})\in D(A^{2})\times D(A)$,
and (\ref{est:GGep'}) in the second case.

\paragraph{\textmd{\textit{Differential inequality for $\EEep$}}}

The time-derivative of (\ref{defn:EEep}) is
\begin{eqnarray}
	\EEep'(t) & = & -\frac{1}{(1+t)^{p}}\frac{|\rep'(t)|^{2}}{\cep(t)}
	\left(2+\ep\frac{\cep'(t)(1+t)^{p}}{\cep(t)}\right)
	\nonumber  \\
	 &  & +\frac{2}{\cep(t)}\langle\rep'(t),\gep(t)\rangle
	+2\langle A\roep(t),\tetep'(t)\rangle.
	\label{comp:EEep}
\end{eqnarray}

By standard inequalities we have that
$$2\langle A\roep(t),\tetep'(t)\rangle\leq
2|A^{1/2}\tetep'(t)|\cdot|A^{1/2}\roep(t)|\leq
2|A^{1/2}\tetep'(t)|\sqrt{\EEep(t)},$$
$$\frac{2}{\cep(t)}\langle\rep'(t),\gep(t)\rangle\leq
\frac{1}{(1+t)^{p}}\frac{|\rep'(t)|^{2}}{\cep(t)}+
\frac{1}{\cep(t)}(1+t)^{p}|\gep(t)|^{2}.$$

Plugging these estimates into (\ref{comp:EEep}), when $\ep$ is small 
enough we obtain that
$$\EEep'(t)\leq 2|A^{1/2}\tetep'(t)|\sqrt{\EEep(t)}+
\frac{1}{\mu}(1+t)^{p}|\gep(t)|^{2},$$
which is exactly (\ref{est:EEep'}) with
$$\psi_{1,\ep}(t):=2|A^{1/2}\tetep'(t)|,
\hspace{3em}
\psi_{2,\ep}(t):=\frac{1}{\mu}(1+t)^{p}|\gep(t)|^{2}.$$

When $(u_{0},u_{1})\in D(A^{3/2})\times D(A^{1/2})$, we have that
$\theta_{0}\in D(A^{1/2})$, hence (\ref{hp:psi12}) follows from
(\ref{th:tetep'-int}) with $k=0$, and (\ref{th:gep-int-1}).

When $(u_{0},u_{1})\in D(A^{2})\times D(A)$, we have that
$\theta_{0}\in D(A)$, and we need (\ref{hp:psi12}) with
$\psi_{1,\ep}(t):=2|A\tetep'(t)|$, and
$\psi_{2,\ep}(t):=\mu^{-1}(1+t)^{p}|A^{1/2}\gep(t)|^{2}$.  Due to the
regularity of $\theta_{0}$, estimate (\ref{hp:psi12}) follows in this
case from (\ref{th:tetep'-int}) with $k=1$, and (\ref{th:gep-int-2}).

\paragraph{\textmd{\textit{Differential inequality for $\FFep$}}}

The time-derivative of (\ref{defn:FFep}) is
\begin{eqnarray}
	\FFep'(t) & = & -\frac{1}{(1+t)^{p}}\frac{|\rep'(t)|^{2}}{\cep(t)}
	\left(2+\ep\frac{\cep'(t)(1+t)^{p}}{\cep(t)}-
	\ep\delta\cep(t)\right)
	\nonumber  \\
	 &  & - \frac{\delta\cep(t)}{(1+t)^{p}}|A^{1/2}\roep(t)|^{2} 
	 -\delta p\frac{|\roep(t)|^{2}}{(1+t)^{2p+1}}
	 -\frac{\ep\delta p}{(1+t)^{1+p}}\langle\rep'(t),\roep(t)\rangle
	\nonumber  \\
	 &  & +\frac{\ep\delta}{(1+t)^{p}}\langle\rep'(t),\tetep'(t)\rangle
	 +2\langle A^{1/2}\roep(t),A^{1/2}\tetep'(t)\rangle
	 +\frac{\delta}{(1+t)^{2p}}\langle\roep(t),\tetep'(t)\rangle
	 \nonumber  \\
	  &  & +\frac{2}{\cep(t)}\langle\rep'(t),\gep(t)\rangle
	  +\frac{\delta}{(1+t)^{p}}\langle\roep(t),\gep(t)\rangle
	  \nonumber \\
	   & = & I_{1}+\ldots+I_{9}.
	  \label{est:I-1-9}
\end{eqnarray}

Let us estimate some of the terms. Clearly we have that $I_{3}\leq 
0$. From (\ref{est:sp-EEep}) we have that
$$I_{6}\leq 2\left|A^{1/2}\roep(t)\right|\cdot
\left|A^{1/2}\tetep'(t)\right|\leq
k_{19}\ep\left|A^{1/2}\tetep'(t)\right|,$$
$$I_{7}\leq\frac{\delta}{(1+t)^{2p}}|\roep(t)|\cdot|\tetep'(t)|\leq
\frac{\delta}{(1+t)^{2p}}\frac{1}{\sqrt{\nu}}
|A^{1/2}\roep(t)|\cdot|\tetep'(t)|\leq k_{20}\ep|\tetep'(t)|.$$

From standard inequalities we have that
$$I_{5}\leq\frac{\ep\delta}{(1+t)^{p}}|\rep'(t)|\cdot|\tetep'(t)|\leq
\frac{\ep\delta}{2}\frac{1}{(1+t)^{p}}\frac{|\rep'(t)|^{2}}{\cep(t)}
+\frac{\ep\delta}{2}\frac{\cep(t)}{(1+t)^{p}}|\tetep'(t)|^{2},$$
$$I_{8}\leq\frac{2}{\cep(t)}|\rep'(t)|\cdot|\gep(t)|\leq
\frac{1}{2}\frac{1}{(1+t)^{p}}\frac{|\rep'(t)|^{2}}{\cep(t)}
+\frac{2}{\cep(t)}(1+t)^{p}|\gep(t)|^{2},$$
$$I_{9}\leq\frac{\delta}{(1+t)^{p}}|\roep(t)|\cdot|\gep(t)|\leq
\frac{\delta\sigma}{2}\frac{1}{(1+t)^{3p}}|\roep(t)|^{2}
+\frac{\delta}{2\sigma}(1+t)^{p}|\gep(t)|^{2}.$$

Plugging all these estimates into (\ref{est:I-1-9}), and recalling
once more assumptions (\ref{hp:cep-below}) through
(\ref{hp:cep'-above}), we obtain that
\begin{eqnarray}
	\FFep'(t) & \leq & -\frac{1}{(1+t)^{p}}\frac{|\rep'(t)|^{2}}{\cep(t)}
	\left(\frac{3}{2}+\ep\frac{\cep'(t)(1+t)^{p}}{\cep(t)}-
	\ep\delta\cep(t)-\frac{\ep\delta}{2}\right)
	\nonumber   \\
	\noalign{\vspace{1ex}}
	 &  & -\frac{\delta\cep(t)}{(1+t)^{p}}|A^{1/2}\roep(t)|^{2} 
	 +\frac{\delta\sigma}{2}\frac{1}{(1+t)^{3p}}|\roep(t)|^{2}
	 -\frac{\ep\delta p}{(1+t)^{1+p}}\langle\rep'(t),\roep(t)\rangle
	 \nonumber   \\
	 \noalign{\vspace{1ex}}
	 &  & +k_{21}\ep\left(|\tetep'(t)|+|\tetep'(t)|^{2}+
	 |A^{1/2}\tetep'(t)|\right)
	 +k_{22}(1+t)^{p}|\gep(t)|^{2}.
	 \label{defn:psi3}
\end{eqnarray}

Let $\psi_{3,\ep}(t)$ denote the sum of the two terms of the last
line.  Then (\ref{est:FFep'}) is proved if we show that the sum of the
terms in the first two lines is less than or equal to
$-\beta(1+t)^{-p}\FFep(t)$ for every $t\geq T$.  In turn, this is
equivalent to showing that $$\hspace{-3em}
\dfrac{|\rep'(t)|^{2}}{\cep(t)}
\left(\frac{3}{2}+\ep\dfrac{\cep'(t)(1+t)^{p}}{\cep(t)}-
\ep\delta\cep(t)-\frac{\ep\delta}{2}-\ep\beta\right)+
(\delta\cep(t)-\beta)|A^{1/2}\roep(t)|^{2}+$$
\begin{equation}
	-\frac{\delta(\sigma+\beta)}{2}\frac{|\roep(t)|^{2}}{(1+t)^{2p}}
	+\left(\dfrac{\ep\delta p}{1+t}
	-\dfrac{\ep\delta\beta}{(1+t)^{p}}\right)
	\langle\rep'(t),\roep(t)\rangle\geq 0
	\label{est:sp-mostro}
\end{equation}
holds true for every $t\geq T$.

Let $S_{1}$, \ldots, $S_{4}$ denote the four terms in 
(\ref{est:sp-mostro}), which we estimate as in the proof of 
Theorem~\ref{thm:lin-h}. From the smallness of $\ep$ we have that
\begin{equation}
	S_{1}\geq\frac{|\rep'(t)|^{2}}{\cep(t)}\geq
	\frac{1}{M_{3}}|\rep'(t)|^{2}.
	\label{est:sp-S1}
\end{equation}

Since $\delta \mu\geq\beta$, from (\ref{hp:coercive}) we have that
\begin{eqnarray*}
	S_{2}+S_{3} & \geq & (\delta \mu-\beta)|A^{1/2}\roep(t)|^{2}-
	\frac{\delta(\sigma+\beta)}{2}\frac{1}{(1+t)^{2p}}|\roep(t)|^{2}
	\\
	 & \geq & \left[(\delta \mu-\beta)\nu
	-\frac{\delta(\sigma+\beta)}{2}\frac{1}{(1+T)^{2p}}\right]
	|\roep(t)|^{2}
\end{eqnarray*}
for every $t\geq T$.  Due to the choices (\ref{choice:sp-p=0}) and
(\ref{choice:sp-p>0}), in both cases the term in brackets is greater
than or equal to $\nu$, hence $S_{2}+S_{3}\geq\nu |\roep(t)|^{2}$ for
every $t\geq T$.  Now we add this inequality to (\ref{est:sp-S1}), and
we apply the inequality between arithmetic and geometric mean, exactly
as in the proof of Theorem~\ref{thm:lin-h}.  If $\ep$ is small enough
we obtain that $$\hspace{-3em} S_{1}+S_{2}+S_{3} \geq
\frac{1}{M_{3}}|\rep'(t)|^{2}+\nu|\roep(t)|^{2} \geq
2\sqrt{\frac{\nu}{M_{3}}}\cdot|\rep'(t)|\cdot|\roep(t)|$$
$$\hspace{3em}
\geq\ep\delta(1+\beta)|\rep'(t)|\cdot|\roep(t)|
\geq \left(\frac{\ep\delta p}{1+t}
+\frac{\ep\delta\beta}{(1+t)^{p}}\right) |\rep'(t)|\cdot|\roep(t)| \geq
|S_{4}|,$$
which proves (\ref{est:sp-mostro}), hence also (\ref{est:FFep'}).

It remains to prove (\ref{hp:psi3}), with $\psi_{3,\ep}(t)$ equal to
the sum of the two terms in the last line of (\ref{defn:psi3}).  

When $(u_{0},u_{1})\in D(A^{3/2})\times D(A^{1/2})$, we have that
$\theta_{0}\in D(A^{1/2})$, hence (\ref{hp:psi3}) follows from
(\ref{th:tetep'-int}) with $k=0$, and (\ref{th:gep-int-1}).

When $(u_{0},u_{1})\in D(A^{2})\times D(A)$, we have that
$\theta_{0}\in D(A)$, and we need (\ref{hp:psi3}) with
$$\psi_{3,\ep}(t):=k_{23}\ep\left(|A^{1/2}\tetep'(t)|+
|A^{1/2}\tetep'(t)|^{2}+|A\tetep'(t)|\right)
+k_{24}(1+t)^{p}|A^{1/2}\gep(t)|^{2}.$$

Due to the regularity of $\theta_{0}$, estimate (\ref{hp:psi3})
follows in this case from (\ref{th:tetep'-int}) with $k=1$, and
(\ref{th:gep-int-2}).

\paragraph{\textmd{\textit{Differential inequality for $\GGep$}}}

The time-derivative of (\ref{defn:GGep}) is
$$\GGep'(t)=-\frac{2}{\ep}\frac{1}{(1+t)^{p}}|\rep'(t)|^{2}
-\frac{2}{\ep}\cep(t)\langle A\roep(t),\rep'(t)\rangle
+\frac{2}{\ep}\langle\gep(t),\rep'(t)\rangle.$$

From standard inequalities we have that
$$-\frac{2}{\ep}\cep(t)\langle A\roep(t),\rep'(t)\rangle\leq
\frac{1}{2\ep}\frac{1}{(1+t)^{p}}|\rep'(t)|^{2}+
\frac{k_{25}}{\ep}(1+t)^{p}|A\roep(t)|^{2},$$
$$\frac{2}{\ep}\langle\gep(t),\rep'(t)\rangle\leq
\frac{1}{2\ep}\frac{1}{(1+t)^{p}}|\rep'(t)|^{2}+
\frac{2}{\ep}(1+t)^{p}|\gep(t)|^{2},$$
hence
$$\GGep'(t)\leq
-\frac{1}{\ep}\frac{1}{(1+t)^{p}}|\rep'(t)|^{2}+
\frac{k_{25}}{\ep}(1+t)^{p}|A\roep(t)|^{2}+
\frac{2}{\ep}(1+t)^{p}|\gep(t)|^{2}.$$

At this point (\ref{est:GGep'}) follows from (\ref{est:sp-order-1})
and (\ref{th:gep-pointwise}).  This completes the proof of
Theorem~\ref{thm:lin-sp}.\qed

\subsection{Proof of Theorem~\ref{thm:lin-opt}}

Let us set
$$\Hep(t):=\left( 
\ep\frac{|\uep'(t)|^{2}}{\cep(t)}+|A^{1/2}\uep(t)|^{2}\right)
\frac{1}{\Phi(t)}
\quad\quad
\forall t\geq 0.$$

Due to (\ref{hp:cep-below}) and (\ref{hp:cep-above}), proving 
(\ref{th:main-opt}) is equivalent to showing that $\Hep(t)\to +\infty$ as 
$t\to +\infty$. Since $(u_{0},u_{1})\neq(0,0)$, the solution is 
nontrivial in the sense that $\Hep(t)>0$ for every $t\geq 0$. Moreover 
we have that
\begin{eqnarray*}
	\Hep'(t) & = & \frac{1}{(1+t)^{p}}\frac{1}{\Phi(t)}
	\frac{\ep|\uep'(t)|^{2}}{\cep(t)}\left(
	-\frac{\Phi'(t)}{\Phi(t)}(1+t)^{p}-
	\frac{2}{\ep}-\frac{\cep'(t)(1+t)^{p}}{\cep(t)}\right) \\
	 &  &
	+\frac{1}{(1+t)^{p}}\frac{1}{\Phi(t)}|A^{1/2}\uep(t)|^{2}\left(
	-\frac{\Phi'(t)}{\Phi(t)}(1+t)^{p}\right).
\end{eqnarray*}

As usual, we have that
$$\frac{|\cep'(t)|(1+t)^{p}}{\cep(t)}\leq\frac{M_{4}}{\mu}.$$

Therefore, assumption (\ref{hp:main-opt}) implies the existence of
$T>0$ (depending on $\ep$, but this is not important) such that
$$-\frac{\Phi'(t)}{\Phi(t)}(1+t)^{p}-
\frac{2}{\ep}-\frac{\cep'(t)(1+t)^{p}}{\cep(t)}\geq 1
\quad\quad\mbox{and}\quad\quad 
-\frac{\Phi'(t)}{\Phi(t)}(1+t)^{p}\geq 1$$
for every $t\geq T$, hence
$$\Hep'(t)\geq\frac{1}{(1+t)^{p}}\Hep(t)
\quad\quad
\forall t\geq T.$$

Since $\Hep(T)>0$, and $p\leq 1$, this differential inequality
implies that $\Hep(t)\to +\infty$ as $t\to +\infty$.\qed

\subsection{Proof of Theorems~\ref{thm:main-p}, \ref{thm:main-h}, 
\ref{thm:main-sp}, \ref{thm:main-opt}}

The existence of solutions to (\ref{pbm:p-eq}), (\ref{pbm:p-data}), 
and (\ref{pbm:h-eq}), (\ref{pbm:h-data}) follows from 
Theorem~\ref{thm:bibl}. Let us set now
$$c(t):=m\left(|A^{1/2}u(t)|^{2}\right),
\hspace{3em}
\cep(t):=m\left(|A^{1/2}\uep(t)|^{2}\right).$$

With a standard approximation procedure, we can assume that
$m(\sigma)$ is of class $C^{1}$, and not just locally Lipschitz
continuous.  As a consequence, also $c(t)$ and $\cep(t)$ are of class
$C^{1}$.  If we show that $c(t)$ and $\cep(t)$ satisfy
(\ref{hp:c-below}) through (\ref{hp:cep-c}), then all conclusions of
Theorems~\ref{thm:main-p} through \ref{thm:main-opt} follow from the
corresponding conclusions of Theorems~\ref{thm:lin-p} through
\ref{thm:lin-opt}.

Assumptions (\ref{hp:c-below}) and (\ref{hp:cep-below}) follow from 
(\ref{hp:ndg}).

Assumptions (\ref{hp:c-above}) and (\ref{hp:cep-above}) follow from
the fact that $|A^{1/2}u(t)|^{2}$ and $|A^{1/2}\uep(t)|^{2}$ are
bounded because of (\ref{th:bibl-p}) and (\ref{th:bibl-h}),
respectively.

Since 
$$c'(t)=2m'\left(|A^{1/2}u(t)|^{2}\right)
\langle Au(t),u'(t)\rangle,$$
assumption (\ref{hp:c'-above}) follows from the boundedness of
$|u'(t)|$, $|A^{1/2}u(t)|$, and $|Au(t)|$, resulting from
(\ref{th:bibl-p}).

Similarly, we have that
$$\cep'(t)=2m'\left(|A^{1/2}\uep(t)|^{2}\right)
\langle A\uep(t),\uep'(t)\rangle,$$
and therefore estimate (\ref{th:bibl-h}) implies that
$$|\cep'(t)|\leq
k_{1}|A\uep(t)|\cdot|\uep'(t)| \leq
k_{2}\frac{1}{(1+t)^{1+p}}\cdot\frac{1}{1+t}\leq
\frac{k_{2}}{(1+t)^{p}},$$
which is exactly (\ref{hp:cep'-above}).

It remains to prove (\ref{hp:cep-c}). To this end, we first remark 
that
\begin{eqnarray*}
	\left|\left|A^{1/2}\uep(t)\right|^{2}-\left|A^{1/2}u(t)\right|^{2}\right|
	& = & \left|\langle A^{1/2}(\uep(t)+u(t)),
	A^{1/2}(\uep(t)-u(t))\rangle\right| \\
	 & \leq & \left(|A^{1/2}\uep(t)|+|A^{1/2}u(t)|\right)
	 \cdot|A^{1/2}\roep(t)|.
\end{eqnarray*}

Now $|A^{1/2}\uep(t)|$ and $|A^{1/2}u(t)|$ are bounded because of
(\ref{th:bibl-p}) and (\ref{th:bibl-h}), and $|A^{1/2}\roep(t)|$ can
be estimated by means of (\ref{th:bibl-sp}).  Since $m(\sigma)$ is
(locally) Lipschitz continuous, we obtain that 
$$|\cep(t)-c(t)|\leq k_{3}
\left|\left|A^{1/2}\uep(t)\right|^{2}-\left|A^{1/2}u(t)\right|^{2}\right|
\leq k_{4}\ep,$$
which is exactly (\ref{hp:cep-c}).\qed

\label{NumeroPagine}

\end{document}